\documentclass[12pt]{article}

\usepackage[margin=1in]{geometry}
\usepackage{hyperref}
\usepackage{amsmath}
\usepackage{amssymb}
\usepackage{amsthm}
\usepackage{amsfonts}
\usepackage{graphicx}
\usepackage{mathrsfs}
\usepackage{color}
\usepackage{accents}
\usepackage{comment}

\newcommand\MR[1]{\relax\ifhmode\unskip\spacefactor3000 \space\fi
	\MRhref{\expandafter\@rst #1 other}{#1}}
\newcommand{\MRhref}[2]{\href{http://www.ams.org/mathscinet-getitem?mr=#1}{MR#2}}

\newtheorem{theorem}{Theorem}[section]

\newtheorem{lemma}[theorem]{Lemma}
\newtheorem{conv}[theorem]{Convention}
\newtheorem{proposition}[theorem]{Proposition}

\newtheorem{definition}[theorem]{Definition}
\newtheorem{remark}[theorem]{Remark}

\newcommand\xor{\mathbin{\triangle}}

\newcommand{\C}{\mathbb{C}}
\newcommand{\D}{\mathbb{D}}
\newcommand{\N}{\mathbb{N}}

\newcommand{\Z}{\mathbb{Z}}
\newcommand{\R}{\mathbb{R}}

\def\bT{\mathbb{T}}
\def\bH{\mathbb{T}^*}

\newcommand{\op}[1]{\operatorname{#1}}
\newcommand{\Ber}{\op{Ber}}
\newcommand{\icon}{\xi_\mesh}
\newcommand{\pivm}{\mu}
\newcommand{\bx}{\mathcal{B}}
\newcommand{\bbH}{\mathbb{H}}
\newcommand{\eps}{\varepsilon}
\renewcommand{\P}{\mathbb{P}}
\newcommand{\E}{\mathbb{E}}

\DeclareMathOperator{\SLE}{SLE}
\DeclareMathOperator{\CLE}{CLE}

\DeclareMathOperator{\Area}{Area} 
\DeclareMathOperator{\Mink}{Mink} 
\DeclareMathOperator{\dist}{dist}
\DeclareMathOperator{\osc}{osc}
\def\om{\mathfrak{m}_\sle} 
\def\omo{\mathfrak{m}_{\sle^0}} 
\def\omp{\mathfrak{m}_{\cP}} 
\def\con1{c_{\op l }}
\def\constp{c_{\op p}}
\def\im{\tau}
\def\mesh{\eta}
\newcommand{\ft}{\mathfrak{t}}
\def\arm{\alpha}
\def \nor{\beta_\eps}
\def \intd{\zeta}

\def\sle{\gamma}
\def\cZ{\mathcal{Z}}

\def\cX{\mathcal{X}}
\def\cW{\mathcal{W}}

\def\cU{\mathcal{U}}

\def\cQ{\mathcal{Q}}
\def\cP{\mathcal{P}}

\def\cN{\mathcal{N}}

\def\cH{\mathcal{H}}
\def\cG{\mathcal{G}}

\def\cD{\mathcal{D}}
\def\cC{\mathcal{C}}
\def\cB{\mathcal{B}}
\def\cA{\mathcal{A}}

\newcommand{\notion}[1]{{\bf  \textit{#1}}}
\newcommand{\1}{\mathbf{1}}
\newcommand{\ol}{\overline}
\newcommand{\wt}{\widetilde}
\newcommand{\ul}{\underline}
\newcommand{\wh}{\widehat}
\newcommand{\ubar}[1]{\underaccent{\bar}{#1}}
\newcommand{\Eb}[1]{\mathbb E\left[{#1}\right]}
\newcommand{\Pb}[1]{\mathbb P\left[{#1}\right]}
\def\H{\mathbb{H}}

\def\Quad{Q}

\def\dist{\mathrm{dist}}
\def \p {{\partial}}
\def\defeq{:=}

\def\mesh{\eta}
\def\A{\mathcal{A}}

\def\@rst #1 #2other{#1}

\newcommand{\arxiv}[1]{\href{http://arxiv.org/abs/#1}{#1}}
\def\MR#1{\href{http://www.ams.org/mathscinet-getitem?mr=#1}{MR#1}}

\begin {document}
\title{Natural parametrization of percolation interface and pivotal points}

\author{
Nina Holden\thanks{
Institute for Theoretical Studies, ETH Z\"urich. Partially supported by a fellowship from the Research Council of Norway and partially supported by Dr.\ Max R\"ossler, the Walter Haefner Foundation, and the ETH Z\"urich Foundation. 
Email: \href{mailto:nina.holden@eth-its.ethz.ch}{nina.holden@eth-its.ethz.ch}.		
}
\and 	
Xinyi Li\thanks{
Beijing International Center for Mathematical Research, Peking University. Research supported by NSFC (No.\ 12071012) and the National Key R\&D Program of China (No.\ 2020YFA0712900).\newline
Email: \href{mailto:xinyili@bicmr.pku.edu.cn}{xinyili@bicmr.pku.edu.cn}.
} 
\and 
Xin Sun\thanks{
Department of Mathematics, University of Pennsylvania. 
Supported by a Junior Fellow award from the Simons Foundation, and NSF Grant DMS-1811092 and DMS-2027986.
Email: \href{mailto:xinsun@sas.upenn.edu}{xinsun@sas.upenn.edu}.
}
}

\date{}
\maketitle

\begin{abstract}
	We prove that the interface of critical site percolation on the triangular lattice converges to $\SLE_6$ in its natural parametrization, where the discrete interface is parametrized such that each edge is crossed in one unit of time, while the limiting curve is parametrized by its $7/4$-dimensional Minkowski content. We also prove that the scaling limit of counting measure on the pivotal points, which was proved to exist by Garban, Pete, and Schramm (2013), is its $3/4$-dimensional Minkowski content up to a deterministic multiplicative constant.

\bigskip
  Nous montrons que l'interface de la percolation du site critique sur le r\'eseau triangulaire converge vers la courbe $\SLE_6$ dans sa parm\'etrisation naturelle, o\`u l'interface discr\`ete est param\'etr\'ee de telle sorte que chaque ar\^ete se croise en une unit\'e de temps, tandis que la courbe limite est param\'etr\'ee par son contenu  $7/4$-dimensionnel de Minkowski. Nous montrons \'egalement que la limite d'\'echelle de la mesure de comptage sur les points pivots, dont l'existence a \'et\'e confirm\'ee par Garban, Pete et Schramm (2013), est son contenu $3/4$-dimensionnel de Minkowski  jusqu'\`a une constante multiplicative d\'eterministe. 
\end{abstract}

\section{Introduction}
Percolation is one of the most studied statistical mechanics models in probability. Since the breakthrough works of Smirnov \cite{Smirnov-Cardy}, who proved the conformal invariance of critical site percolation on the triangular lattice, and of Schramm \cite{Schramm-SLE}, who introduced the Schramm-Loewner evolution (SLE), the understanding of the scaling limit of percolation on planar lattices has greatly improved. 

Garban, Pete, and Schramm \cite{GPS-piv,gps-fourier,gps-near-crit} made important contributions in this direction. In \cite{GPS-piv} they proved scaling limit results for several important classes of points for critical percolation, including pivotal points and points on the percolation interface. They proved that the limiting measures are conformally covariant, and that they are measurable with respect to the scaling limit of percolation.

In the continuum, a substantial effort has been made to understand natural measures on special points of SLE$_\kappa$ curves. For example, SLE$_\kappa$ curves 
have non-trivial
$2\wedge (1+\kappa/8)$-dimensional Minkowski content, which defines a parametrization of the curve called the \notion{natural parametrization}.
SLE with its natural parametrization is uniquely characterized by conformal invariance and domain Markov property, with the constraint that the parametrization is rescaled in a covariant way under the application of a conformal map. See \cite{Lawler-Sheffield,Lawler-Zhou,Lawler-Minkowski}.
SLE with its natural parametrization is believed to describe the scaling limit of curves in statistical physics models parametrized such that one edge/face/vertex  is visited in one unit of time. This conjecture was proved for the case of the loop-erased random walk (LERW) and SLE$_2$ by Lawler and Viklund \cite{lawler-viklund-lerw-nat,lawler-viklund-lerw-radial}.

In this paper, we link the limiting measures in \cite{GPS-piv} with the natural measures on special points of $\SLE_6$. The purpose of building the link is two-fold:
\begin{enumerate}
	\item It makes the limiting measures in \cite{GPS-piv} more intrinsic and concrete. In the case of percolation pivotal points, this link is important for the work of the first and third authors on the conformal embedding of uniform triangulations~\cite{Cardy}.
	\item It allows us to prove that the percolation interface converges to $\SLE_6$ in its natural parametrization. 
\end{enumerate}

\subsection{The scaling limit of the percolation interface under its natural parametrization}
\label{sec:nat}

Let us briefly recall the definition of SLE. Fix $\kappa>0$, and let $(B_t)_{t\geq 0}$  be a standard linear Brownian motion. Consider the Loewner differential equation 
\[
\partial_t g_t(z) =\frac{2}{g_t(z)-\sqrt \kappa B_t},\qquad g_0(z)=z, \; \forall z\in \ol\bbH.
\] Then for each $z\in\H$, $g_t(z)$ is well-defined up to some time $\tau_z\in[0,\infty]$. Let $K_t=\ol{\{z:\tau_z<t\}}$. Then a.s.\ there exists a unique continuous non-self-crossing curve $\gamma$ such that
$K_t$ is the closure of points disconnected from $\infty$ on $\bbH$ by $\gamma([0,t])$. 
We call $\gamma$ the chordal $\SLE_\kappa$ on $\mathbb H$ from $0$ to $\infty$ (under the capacity parametrization).
Let  $\Omega$ be a simply connected domain whose set of prime ends $\p \Omega$ is a continuous image of a circle\footnote{This is the necessary and sufficient  boundary condition for the Riemann mapping from the unit disk to $\Omega$ to continuously extend to the boundary. (See e.g. \cite{Riemann-Boundary}.)}. Let $a,b$ be two distinct points on $\p \Omega$.
Consider  a conformal map $f:\bbH\to \Omega$ with $f(0)=a$ and $f(\infty)=b$. Although there is one degree of freedom when choosing $f$, the law of $f(\sle)$ (viewed as a continuous curve modulo increasing reparametrizations)  does not depend on this choice. We call this probability measure the \notion{chordal $\SLE_\kappa$ on $\Omega$ from $a$ to $b$}, or simply $\SLE_\kappa$ on $(\Omega,a,b)$.

Let $\bT$ denote the regular  triangular lattice where each face is an equilateral triangle. 
For $\mesh>0$, let $\mesh \bT$  be $\bT$  rescaled by $\mesh$. Each vertex on $\mesh\bT$ is called a site.
Let $\mesh\bH$ denote the regular hexagonal lattice dual to $\mesh\bT$ such that each vertex on $\bT$ corresponds to a hexagonal face on $\mesh\bH$.
Given an edge $e$ of $\mesh\bT$, let $e^*$ be its dual edge in $\mesh\bH$. Recall that a Jordan domain is a bounded simply connected domain on $\C$ whose boundary is homeomorphic to a circle. A Jordan domain $D$ is called a \notion{$\mesh$-polygon} if $\p D$ lies on the lattice $\mesh\bT$. A vertex $v$ on $\mesh\bT$ is called an \notion{inner vertex} (resp., \notion{boundary vertex}) of $D$ if $v\in D$ (resp., $v\in\p D$). We similarly define \notion{boundary/inner edges} of $D$. 

Suppose $\Omega$ is  a Jordan domain.  Let $\Omega_\mesh$ be the largest $\mesh$-polygon whose set of inner vertices is contained in $\Omega$ and forms a connected set on $\mesh\bT$. (In case of a draw, choose $\Omega_\mesh$ arbitrarily from the set of largest $\mesh$-polygons, but note that $\Omega_\mesh$ will be uniquely determined for all sufficiently small $\mesh$.)    Including all inner vertices and edges of $\Omega_\mesh$, we obtain a planar graph embedded in $\C$ which we  call the \notion{$\mesh$-approximation} of $\Omega$ and still denote by $\Omega_\mesh$.%
\footnote{A notion of $\mesh$-approximation of the Jordan domain $\Omega$ is also introduced  in Definition~4.1 of \cite{Camia-Newman-CLE1}, which is denoted by $D^\mesh$ in their notation.
One can check that  $D_\mesh$  equals the union of $D^\mesh$ and its so-called external boundary  defined in \cite[Section 4]{Camia-Newman-CLE1}.}
To distinguish with the continuum, we write the union of boundary vertices and edges of $\Omega_\mesh$ as $\Delta\Omega_\mesh$.  A \notion{path} on a graph is a sequence of vertices such that each vertex is adjacent to its successor. 
Given two distinct boundary edges of  $\Omega_\mesh$, removing $\{e,e'\}$ from $\Delta\Omega_\mesh$ gives two paths on the boundary.
We let $\Delta_{e,e'}\Omega_\mesh$ denote the one tracing $\Delta\Omega_\mesh$ counterclockwise from $e$ to $e'$.
Given  $x\in \p\Omega$, let $x_\mesh$ be the  edge on $ \Delta\Omega_\mesh$ closest to $x$ (if there is a tie, choose one arbitrarily). 

A \notion{site percolation} on $\Omega_\mesh$ is a  black/white coloring of inner vertices of $\Omega_\mesh$.
The \notion{critical  Bernoulli site percolation} on $\Omega_\mesh$, which we denote by $\Ber(\Omega_\mesh)$,  is the uniform measure on site percolations on $\Omega_\mesh$.
A  coloring of vertices on $\Delta\Omega_\mesh$ is called a \textit{boundary condition}.  
A site percolation on $\Omega_\mesh$ together with a  boundary condition determines  a coloring of vertices on $\Omega_\mesh$. 
The \notion{$(a,b)$-boundary condition} is the coloring where vertices on $\Delta_{a_\mesh,b_\mesh}\Omega_\mesh$ (resp., $\Delta_{b_\mesh,a_\mesh}\Omega_\mesh$) are black (resp., white).  Note that this is well-defined since we required that $\Omega_\eta$ does not have any cut vertices.
Given a site percolation  $\omega_\mesh$ on $\Omega_\mesh$ with $(a,b)$-boundary condition,
there is a unique  path  $\sle_\mesh$ on $\mesh\bH$ from $a^*_\mesh$ to $b^*_\mesh$, such that each edge on the path has a white vertex on its left side and a black vertex on its right side.
We call $\sle_\mesh$ the \notion{percolation  interface} of $\omega_\mesh$  on  $(\Omega_\mesh,a_\mesh,b_\mesh)$.

Let $(\cU,d_\cU)$ denote the separable metric space of continuous curves modulo reparametrization, with the distance $d_\cU$ between curves $\gamma^1:[0,T_1]\to\C$ and 
$\gamma^2:[0,T_2]\to\C$ defined by
\begin{equation}\label{eq:u-metric}
	d_\cU (\gamma_1,\gamma_2) = \inf_{\alpha,\beta} \left[ \sup_{0\le t\le 1}\big|\sle^1(\alpha(t))  -\sle^2(\beta(t))\big|  \right],
\end{equation}
where the infimum is taken over all choices of increasing bijections $\alpha:[0,1]\to [0,T_1]$ and 
$\beta:[0,1]\to [0,T_2]$. It is  proved in \cite{Smirnov-Cardy,CN-SLE} that $\sle_\mesh$ converges to an $\SLE_6$ on $(\Omega,a,b)$ for the $d_\cU$-metric (see Theorem~\ref{thm:embedding}). Although this convergence result gives a powerful tool for analyzing large scale properties of percolation (e.g.\ arm exponents \cite{Sminov-Werner}), a more natural notion of convergence would be under the parametrization where $\sle_\mesh$ traverses each edge in the same amount of time. We prove this result in Theorem \ref{thm:conv} below. Before stating this result, we need  the notions of Minkowski content and occupation measure.

\begin{definition}\label{def:Mink}
	Given a set $\cA\subset \C$, for $r>0$, let $\cA^r=\{z\in\C: B(z,r)\cap \cA\neq \emptyset \}$. For $d\in [0,2]$ we define the  \notion{$d$-dimensional Minkowski content} of $\cA$ to be the following limit, provided it exists
	\begin{equation}\label{eq:Mink}
	\Mink_d(\cA):=\lim_{r\to 0} r^{d-2} \Area(\cA^r).
	\end{equation} 
	If the limit does not exist, then the $d$-dimensional Minkowski content of $\cA$ is not defined.
\end{definition}

\begin{definition}\label{def:occu}
Fix $d\in [0,2]$.  Let $\cA\subset \C$ be a random closed set and let $\mu$ be a random Borel measure on $\C$. 
Suppose $\P[\mu(U)=\Mink_d(\cA\cap U)]=1$ for each Jordan domain $U$ with piecewise smooth boundary. We call  $\mu_\cA$ the \notion{occupation measure} of $\cA$ and say that it is (a.s.) defined by the $d$-dimensional Minkowski content of $\cA$.
\end{definition}
Let $\gamma$ be an $\SLE_6$ on $(\Omega,a,b)$, where $\Omega$ is a Jordan domain with smooth boundary and $a,b$ are two distinct boundary points. 
Assume the parametrization of $\gamma$ comes from the image of a capacity-parametrized $\SLE_6$ on $(\bbH, 0,\infty)$ under a conformal map $f:\bbH\to \Omega$ with $f(0)=a,f(\infty)=b$. By \cite{Lawler-Minkowski}, we know the following.
\begin{enumerate}
	\item\label{item:occupation} For each $t\in(0,\infty]$, a.s.\ the $7/4$-dimensional Minkowski content of $\sle([0,t])$ exists and defines the occupation  measure of  $\gamma([0,t])$ as in Definition~\ref{def:occu}. We denote the occupation measure of $\gamma((0,\infty))$ by $\om$.
	\item\label{item:regular}
	The function $t\mapsto\Mink_{7/4}\left(\sle([0,t])\right)$ is a.s.\  strictly increasing and H\"{o}lder continuous.
\end{enumerate} 
\begin{definition}\label{def:natural}
Suppose $\Omega$ is a Jordan domain with smooth boundary. Let $\gamma$ be an $\SLE_6$ on $(\Omega,a,b)$, where $a,b$ are two distinct boundary points.
Let $\wh\sle:[0, \om(\Omega)] \to\Omega$ be the parametrization of $\sle$ such that $\Mink_{7/4}(\sle([0,t]))$ $=t$ for any $t\in [0, \om(\Omega)]$. 
Then $\wh\sle$ is called the \notion{natural parametrization} of $\sle$. 
\end{definition}
In fact \cite{Lawler-Minkowski} mainly focuses on  the upper half plane. However, as explained below Theorem~1.1 there, the case of Jordan domains with smooth boundary can be easily obtained by the covariance of Minkowski content under conformal mappings. 

Define the following distance $\rho$ between two parametrized curves $\sle^1:[0,T_1]\to \C$ and $\sle^2:[0,T_2]\to \C$.
\begin{equation}\label{eq:distance}
	\rho(\sle^1,\sle^2) = \left[|T_2-T_1|+ \sup_{0\le s\le 1}|\sle^2(s T_1)  -\sle^1(s T_2)|  \right].
\end{equation}
As mentioned above, for statistical mechanical models where SLE is the scaling limit in the $d_\cU$-metric, 
it is believed that the convergence should also hold in the $\rho$-metric under the natural parametrization.
In this paper, we prove this for the percolation interface.
\begin{theorem}\label{thm:conv} 
Let $(\Omega,a,b)$, $\sle$, and $\wh\sle$ be as in Definition~\ref{def:natural}. 
For $\eta>0$, sample $\omega_\mesh$ from $\Ber(\Omega_\mesh)$ and let
$\sle_\mesh$ be the  interface of $\omega_\eta$ from $a_\mesh$ to $b_\mesh$.
{Pick $\con1>0$, write $\icon=\con1 \eta^2/\alpha_2^\eta(\eta,1)$, and} let 
$\wh \sle_\mesh$ be the parametrization of $\sle_\mesh$ with constant speed (with respect to the Euclidean metric) such that each edge is crossed in $\icon$ units of time. Then with an appropriate choice of $\con1$, the curve $\wh \sle_\mesh$ converges weakly to $\wh\sle$  in the $\rho$-metric.
\end{theorem}
Fix $\con1>0$, and let the (normalized) interface measure  $\im_\mesh$ on $\sle_\mesh$ be defined by 
\begin{equation}\label{eq:im}
\im_\mesh \defeq \con1\sum_{e\in \sle_\mesh} \delta_e \frac{\mesh^2}{\arm^\mesh_2(\mesh,1)},
\end{equation}
where $\arm^\mesh_2(\mesh,1)$ is a normalizing constant depending on $\mesh$ that we will specify in Section~\ref{subsec:arm}, and 
$\delta_e$ is the measure assigning unit mass uniformly along $e$
and 0 elsewhere. The following is proved  in \cite{GPS-piv}.
\begin{theorem}[\cite{GPS-piv}]\label{thm:length}
In the  setting of Theorem~\ref{thm:conv}, there is a coupling of   $(\omega_\mesh)_{\mesh>0}$ and $\gamma$ such that as $\mesh\to 0$, it holds a.s.\  that $\sle_\mesh$ converges to $\gamma$
in the $d_\cU$-metric, and $\im_\mesh$ in \eqref{eq:im} converges to a random Borel measure  $\im$  supported on the range of $\sle$ in the weak topology. Moreover, $\tau$ is measurable with respect to $\sle$.
\end{theorem}
It was not proved in \cite{GPS-piv} that the measure $\tau$ defines a parametrization of $\gamma$. 
Some of the challenges in proving this are discussed in \cite[Sections 1.2 and 5.3]{GPS-piv}.

As a first step towards proving Theorem~\ref{thm:conv}, we prove the following in Section~\ref{sec:equiv}.
\begin{theorem}\label{thm:equi-length}
In Theorem~\ref{thm:conv}, one can choose $\con1$  in \eqref{eq:im} such that $\im=\om$ a.s.
\end{theorem}

{We end this subsection by commenting on our proof ideas for Theorems~\ref{thm:equi-length} and~\ref{thm:conv}. 
On the one hand, the proof of  Theorem~\ref{thm:equi-length}  closely follows~\cite[Section~4]{GPS-piv} with the simplification that we work directly in the continuum, hence our one-point and two-point estimates are  power laws with no sub-polynomial corrections, in contrast to the arm exponent estimates for percolation. 
On the other hand, an additional technicality arises when we try to implement the continuum analog of a strong coupling result from~\cite[Section~4]{GPS-piv}. 
See the beginning of Section~\ref{sec:equiv} for more discussion. Our proof of Theorem~\ref{thm:piv} below uses the same idea as in Theorem~\ref{thm:equi-length}. 
However, Theorem~\ref{thm:equi-length} itself is not sufficient for proving Theorem~\ref{thm:conv}  due to the presence of double points in $\SLE_6$. 
To  deal with this issue, we prove that the occupation measure of the frontier of $\SLE_6$ is 0 and use it to conclude the proof of Theorem~\ref{thm:conv}  in Section~\ref{sec:conv}.
}

\subsection{The natural measure on pivotal points}\label{subsec:piv}
Let $\Omega$ be a Jordan domain with smooth boundary and sample $\omega_\mesh$ from $\Ber(\Omega_\mesh)$.
Let $a,b,c,d\in \p \Omega$ be four distinct points ordered counterclockwise. For $\mesh$ small enough such that $a_\mesh$, $b_\mesh$, $c_\mesh$, and $d_\mesh$ are distinct.
the following three sentences describe the same event.
\begin{itemize}
\item  There is a path $\{v_i\}_{1\le i\le n}$ such that  $v_1$ and $v_n$ are on $\Delta_{b_\mesh, c_\mesh}\Omega_\mesh$ and  $\Delta_{d_\mesh, a_\mesh}\Omega_\mesh$ respectively, while $v_i$ is a white inner vertex for all $1<i<n$.
\item Let $e_\mesh$ be the first edge crossed by the percolation interface on $(\Omega_\mesh, a_\mesh, c_\mesh)$  with one endpoint lying on $\Delta_{b_\mesh,d_\mesh} \Omega_\mesh$. Then $e_\mesh$ has an endpoint on   $\Delta_{b_\mesh,c_\mesh} \Omega_\mesh$.
\item 
Let $\ol e_\mesh$ be the first edge crossed by the percolation interface on $(\Omega_\mesh, c_\mesh, a_\mesh)$  with one endpoint lying on $\Delta_{d_\mesh,b_\mesh} \Omega_\mesh$. Then $\ol e_\mesh$ has an endpoint on   $\Delta_{d_\mesh,a_\mesh} \Omega_\mesh$.
\end{itemize}
Denote this event by $E_\mesh$.
Consider the pair of curves $(\gamma^1_\mesh,\gamma^2_\mesh)$ defined as follows. 
When $E_\mesh$ occurs, let $\gamma^1_\mesh$ and $\gamma^2_\mesh$ be the percolation interfaces on $(\Omega, a_\mesh,b_\mesh)$ and $(\Omega_\mesh, c_\mesh,d_\mesh)$, respectively. Otherwise, let $\gamma^1_\mesh$ and $\gamma^2_\mesh$ be the percolation interfaces on $(\Omega, a_\mesh,d_\mesh)$ and $(\Omega_\mesh, c_\mesh,b_\mesh)$, respectively.
Given an event defined in terms of $\omega_\eta$, a site in $\Omega_\mesh$ is called a \notion{pivotal point} for this event   if flipping the color of the site changes the outcome of the event.  Let $\cP_\mesh$ be the set of pivotal  points  for $E_\mesh$. Then a site of $\Omega_{\mesh}$ belongs to $\cP_\mesh$ if and only if it is the endpoint of one edge crossed by $\gamma^1_\mesh$ and one edge crossed by $\gamma^2_\mesh$.

The picture above has a natural scaling limit. Let $\p_{a,b}\Omega$ be the counterclockwise arc on $\p \Omega$ between $a$ and $b$. By locality, we can couple the chordal $\SLE_6$ on $(\Omega,a,b)$ to the chordal SLE$_6$ on $(\Omega,a,d)$ such that the two curves agree  until hitting  the arc $\p_{b,d}\Omega$, after which they evolve independently.  Let  $E$ be the event that the hitting location on $\p_{b,d}\Omega$  lies  on $\p_{b,c}\Omega$. If $E$ occurs (resp., does not occur), let $\gamma^1$ be the $\SLE_6$ from $a$ to $b$ (resp., $d$) so that there exists a unique connected component of $\Omega\setminus  \gamma^1$ whose boundary contains $c$ and $d$ (resp.\ $b$). 
Conditioning on $\gamma^1$, 
let $\gamma^2$ be a chordal $\SLE_6$ on this component from $c$ to $d$ (resp.\ $b$).

A point is called a \notion{pivotal point} for $E$ if and only if it is on the range of both $\sle^1$ and $\sle^2$. Let $\cP$ denote the set of pivotal points of $E$.
Fix $\constp>0$ and define 
\begin{equation}\label{eq:pm} 
	\pivm_\mesh \defeq \constp\sum_{z\in \cP_\mesh} \delta_z \frac{\mesh^2}{\arm^\mesh_4(\mesh,1)},
\end{equation}
where $\arm^\mesh_4(\mesh,1) $   is a normalizing constant which will be specified in Section~\ref{subsec:arm}. 

The following theorem follows from  \cite{GPS-piv}. 
\begin{theorem}\label{thm:length2}
	There is a coupling of $(\gamma^1_\mesh,\gamma^2_\mesh)$ and $(\gamma^1,\gamma^2)$ such that a.s.,
	\begin{itemize}
		\item[(1)] $\1_{\E_\mesh}$ converges to $\1_{E}$,
		\item[(2)] $(\gamma^1_\mesh,\gamma^2_\mesh)$ converges to $(\gamma^1, \gamma^2)$ in the $d_\cU$-metric, and 
		\item[(3)] $\pivm_\eta$ converges to  a measure $\pivm$ supported on $\cP$.
	\end{itemize}
\end{theorem}
Certain basic properties of $\mu$ were also proved in \cite{GPS-piv}, for example that $\mu$ is measurable with respect to the scaling limit of percolation in quad-crossing space (see Section \ref{sec:pre} for definitions) and the conformal covariance of $\mu$.

{As we will explain in more detail in Section~\ref{sec:piv}, (see the discussion above Lemmas~\ref{lem:cut} and~\ref{lem:coupling})}
the set $\cP$ is locally absolutely continuous with respect to the set of cut points of two-dimensional Brownian motion, whose occupation measure is the subject of \cite{cutpoint}. 
{Using the relationship between $\cP$ and Brownian cut points we will prove the following result in Section~\ref{sec:piv}.}
\begin{proposition}\label{prop:occu-piv}
The occupation measure $\omp$ of $\cP$ a.s.\ exists and is defined by its $3/4$-dimensional Minkowski content in the sense of Definition~\ref{def:occu}.
\end{proposition}
In Section~\ref{sec:piv}, we use the same arguments as for Theorem~\ref{thm:equi-length} to conclude the following.
\begin{theorem}\label{thm:piv}In Theorem~\ref{thm:length2}, one can choose $\constp$  in \eqref{eq:pm} such that $\mu=\omp$ a.s.
\end{theorem}
Theorem~\ref{thm:piv} confirms that the scaling limit of the  pivotal measure in \cite{GPS-piv} is in fact induced by the $3/4$-dimensional Minkowski content of the continuum pivotal points. We can also consider double points of $\SLE_6$ and the points of intersection of $\CLE_6$ loops, which describe the full scaling limit of the interfaces between black and white clusters in critical percolation \cite{Camia-Newman-CLE1}.  
In these cases, the analog of Theorem~\ref{thm:piv} holds since their local pictures are absolutely continuous with respect to each other and the Minkowski content is defined locally. We restrict to the formulation in Theorem~\ref{thm:piv} for concreteness.

{Theorem~\ref{thm:piv} is an important ingredient of the first and third authors' proof~\cite{Cardy} of the convergence of uniform triangulations to Liouville quantum gravity (LQG) with parameter $\sqrt{8/3}$ under the so-called Cardy embedding, which is a discrete conformal embedding based on percolation. A key tool in the proof is the 
	\emph{Liouville dynamical percolation} (LDP) introduced in~\cite{ghss-LDP}, which is a variant of  the ordinary dynamical percolation considered in~\cite{gps-near-crit}. The discrete (ordinary) dynamical percolation in~\cite{gps-near-crit}  is defined as follows. We start from a sample of critical Bernoulli site percolation and then use i.i.d.\ exponential clocks at each site to update the color. The discrete LDP is defined in the same way except that the rates of the exponential clocks are not identical but depend on {a background LQG surface}. The continuous LDP is the continuum limit of the discrete LDP as the lattice size and the clock rates are rescaled {appropriately}.  By the existence of the scaling limit of the pivotal measure from~\cite{GPS-piv}, the existence of the continuous LDP was proved in~\cite{ghss-LDP} in the so-called quad-crossing topology (see Section~\ref{subsec:quad}), similarly as in~\cite{gps-near-crit}.

The key idea of~\cite{Cardy} is to consider ordinary dynamical percolation (namely, with i.i.d.\ clocks) on a uniform triangulation and realize that under the conformal embedding the scaling limit of this dynamic is the continuous LDP. Once this is proved, ergodicity of continuous LDP proved in~\cite{ghss-LDP} implies that uniform triangulations under the Cardy embedding converge to LQG. The quad-crossing topology allows~\cite{ghss-LDP} to apply the powerful machinery of noise sensitivity developed in~\cite{gps-fourier} to prove the desired ergodicity of continuous  LDP. However, it is not the natural topology to describe the scaling limit of the ordinary dynamical percolation on uniform triangulations. The natural topology is given by the  \emph{mating-of-trees} framework of Duplantier, Miller, and Sheffield~\cite{dms-mating}; also see~\cite{ghs-survey}.

{The technical bulk of~\cite{Cardy} is to show that the quad-crossing and mating-of-trees descriptions of continuous LDP are equivalent.
In both descriptions, the dynamic is determined by its initial configuration and a Poisson point process whose intensity measure {is supported on the set of pivotal points}, and the equivalence of the LDP descriptions can therefore be reduced to the equivalence of two notions of pivotal measure. The notion of pivotal measure coming from~\cite{ghss-LDP} is given by the ordinary pivotal measure in~\cite{GPS-piv} weighted by the exponential of a Gaussian free field. The notion coming from mating-of-trees, which is introduced in~\cite{bhs-site}, is defined using Brownian motion and involves neither the ordinary pivotal measure from~\cite{GPS-piv} nor Gaussian free field. To show the equivalence of these two notions of pivotal measure, the description of ordinary pivotal measure in terms of Minkowski content as in Theorem~\ref{thm:piv} plays an important role. In particular, with this definition the equivalence of the two pivotal measures becomes a natural and concrete statement for $\CLE_6$.
See~\cite[Section~5]{Cardy} for the detail of this argument and see~\cite[Section~1.4]{Cardy} for an overview of the entire program.}

}

\section{Preliminaries}\label{sec:pre}
In this section we review some basic facts about percolation which are used in later proofs. Most facts are either known or easy consequences of known results. Therefore we will be brief and  refer to \cite{Sminov-Werner,Smirnov-Schramm,GPS-piv,Werner-Percolation-Notes} for more details. 
\subsection{Basic notations}
Throughout the paper, we use $\sle$ and $\sle_\mesh$ to represent $\SLE_6$ and the percolation interface, respectively.  Both $\sle$ and $\sle_\mesh$ are understood as continuous curves modulo reparametrization unless otherwise specified. When there is no risk of confusion, we also use $\sle,\sle_\mesh$ to denote the range of the curves. 

For all $R>0$ and $z\in \C=\R^2$, we let $\bx_R(z)=z+[-R,R]^2$ denote the square of side length $2R$ centered at $z$. We call a set a \notion{box} if it can be written on this form. We write $\bx_R$ for $\bx_R(0)$ and $c\bx_R(z)$ for $z+[-Rc,Rc]^2$
	(instead of $cz+[-Rc,Rc]^2$)
	 for all $c>0$.
For $0<r<R$, let $A(r,R)= \bx_R\setminus{\bx_r}$. We call a domain $A$ an \notion{annulus} if $A$ is topologically equivalent to $A(1,2)$, and we  use $\p_1 A$ and $\p_2 A$ to denote its inner and outer boundaries, respectively. 

Given any two sets $X,Y\subset \R^2$, we write $\dist(X,Y)\defeq\inf\{|x-y|: x\in X,\, y\in Y\}$. Let $\ol X$ denote the closure of $X$. If  $\ol X\subset Y$, we write $X\Subset Y$.

We use classical asymptotic notations. Given  two non-negative functions $f$ and $g$, we write $f\lesssim g $ (resp., $f\gtrsim g $) if there is a constant $C>0$ such that $f(x)\le  Cg(x)$ (resp., $f(x)\ge  Cg(x)$) for all $x$. We also write $f=O(g)$ when $f\lesssim g$.
We write $f\asymp g$ if $f\lesssim g$ and $g\lesssim f$.
We say $f(x)=o_x(1)$ as $x\to a$ if $\lim_{x\to a} f(x)=0$.

\subsection{Quad-crossing representations of percolation}\label{subsec:quad}
There are various ways to represent the scaling limit of critical planar percolation (see e.g.\ the introduction of \cite{Smirnov-Schramm}).
One way is to use its crossing information, as we review now.

A \notion{quad} in $\C$ is
a homeomorphism $Q:[0,1]^2\to\C$.  Let
\begin{align*}
	\partial_1 Q\defeq Q(\{0\}\times [0,1]), \quad\partial_2 Q\defeq Q( [0,1]\times \{0\}),\\
	\partial_3 Q\defeq Q(\{1\}\times [0,1]), \quad \partial_4 Q\defeq Q([0,1]\times \{1\}).
\end{align*}
We will identify a quad $Q$ with $(Q[0,1]^2, Q(0,0), Q(1,0) , Q(1,1), (0,1))$, so quads giving the same such tuple are identified.
Let $\cQ$ be the space of quads in $\C$, equipped with the uniform topology. 
A \notion{crossing} of a quad $Q$ is a closed set in $\C$ containing a connected closed subset of $Q([0,1]^2)$ that intersects both $\partial_1 Q$ and $\partial_3 Q$. 
Given $Q_1,Q_2$ in $\C$, we say $Q_1 \le Q_2$ if every crossing of $Q_2$ contains a crossing of $Q_1$. We say $Q_1<Q_2$ if there exists a neighborhood of $\cN_i$ ($i=1,2$) of $Q_i$ in $\cQ$ such that $N_1\le N_2$ for any $N_i\in \cN_i$.
A \notion{quad-crossing configuration} on $\C$ is a function $\omega:\cQ\to \{0,1\}$ such that the set $\omega^{-1}(1)$ is  closed in $\cQ$ and  for any $Q_1,Q_2$ with  $Q_1< Q_2$, we have  $\omega(Q_2)\le \omega(Q_1)$. We denote the space of quad-crossing configurations on $\C$ by $\cH$.
The set $\cH$ can be endowed with a metric $d_\cH$ such that $(\cH,d_\cH)$  is compact and separable.

Let $\Omega\subsetneqq\C$ be an open set and let $\cQ_\Omega$ be the space of quads with image in $\Omega$. By restricting to $\cQ_\Omega$, each element in $\cH$ induces a quad-crossing configuration on $\Omega$. Let $\cH_\Omega$ be the space of such configurations, endowed with the metric induced by $d_\cH$, which we still denote by $d_\cH$.  We refer to \cite{Smirnov-Schramm,GPS-piv} for more details on $(\cH_\Omega, d_\cH)$. Here we only record the following facts. Suppose $\Omega$ is a Jordan domain and that $\omega_\mesh$ is sampled from $\Ber(\Omega_\mesh)$.  We identify $\omega_\mesh$ with  an element in $\cH_\Omega$ by setting $\omega_\mesh(Q)=1$ if and only if the white sites of $\omega_\mesh$ form a crossing of  $Q$.  Then $\omega_\mesh$ weakly converges to a random variable $\omega$  in $\cH_\Omega$ under the $d_\cH$-metric.  Moreover, 
\begin{enumerate}
	\item for each deterministic quad $Q\in \cQ_\Omega$, in any coupling where $\omega_\mesh\to \omega$ a.s., we have $\omega_\mesh(Q)\to \omega(Q)$ in probability;
	\item there exists a countable collection $\{Q_n\}_{n\in \N}\subset \cQ_\Omega$ such  that $Q_n$ has piecewise smooth boundary and $\{\omega(Q_n)\}_{n\in \N}$ generates the Borel $\sigma$-field of $(\cH_\Omega,d_\cH)$.
\end{enumerate}

\subsection{Some scaling limit results}\label{sub:SLE}
The following scaling limit result is from \cite{Camia-Newman-CLE1} and \cite{GPS-piv}. 
\begin{theorem}\label{thm:embedding}
	Suppose $\Omega$ is a Jordan domain. Sample $\omega_\mesh$ from $\Ber(\Omega_\mesh)$. 
	Then there is a coupling of $(\omega_\mesh)_{\mesh>0}$ such that the following hold.
	\begin{enumerate}
		\item  For any fixed $x,y\in \partial \Omega$ with $x\neq y$,  the interface $\sle^{xy}_\mesh$ on $(\Omega_\mesh,x_\mesh,y_\mesh)$ converges in probability to an $\SLE_6$ curve $\sle^{xy}$ on $(\Omega,x,y)$ under the $d_{\cU}$-metric.
		\item The quad-crossing configuration  $\omega_\mesh$ converges to $\omega$ in probability under the $d_\cH$-metric. 
	\end{enumerate}
In particular,	this provides a coupling of $\omega$ and $\{\gamma^{xy}: x\neq y,\;x,y\in\p\Omega\}$.
\end{theorem}
Theorem~\ref{thm:embedding} is obtained by considering  the collection of disjoint loops $\Gamma_\mesh$ which are interfaces between black and white clusters of $\omega_\mesh$. They converge to
a random collection of loops $\Gamma$ called the conformal loop ensemble with $\kappa=6$ (CLE$_6$) on $\Omega$. Moreover, both $\{\sle^{xy} \}_{x,y\in\p\Omega}$ and $\omega$ are measurable with respect to  the $\CLE_6$ (see \cite{Camia-Newman-CLE1} and \cite[Section~2.3]{GPS-piv}).
We will not give more detail on $\CLE_6$ as it is not needed, but refer to \cite{Camia-Newman-CLE1,Sheffield-CLE} for further details. 

We have seen several  classes of domains so far. For the definition of $\SLE_6$, we assumed that the boundary is a continuous image of a circle. For quad-crossing space, we considered general domains. For Theorem~\ref{thm:embedding}, we considered Jordan domains. In Theorems~\ref{thm:conv} and~\ref{thm:piv},   we assumed that $\p \Omega$ is a smooth Jordan curve.
We  will carefully organize the argument so that  Theorem~\ref{thm:embedding} does not have to be extended to domains with  rougher boundary. See Remark~\ref{rmk:boundary}.

The following gives the convergence of the interface at the hitting time of certain domains. {The lemma will be used to prove Lemmas \ref{lem:GG}, \ref{lem:GG2}, and \ref{lem:tight}.}
\begin{lemma}\label{lem:stop}
	In the setting of Theorem~\ref{thm:conv},  view $\sle,\sle_\mesh:[0,1]\to\Omega$ as parametrized curves coupled together  such that  $\lim_{\mesh\to 0} \sup_{0\le t\le 1}\{ |\sle_\mesh(t) -\sle(t)|\}=0$ a.s. (The existence of such  parametrizations  and couplings  is guaranteed by $d_\cU$-convergence of $\sle_\mesh$ to $\sle$.)
	Let $\sigma_\mesh,\sigma$ be stopping times for $\sle_\mesh$ and $\sle$, respectively, such that   $\sigma_\mesh\to\sigma$ a.s. 
	Fix a piecewise smooth simple curve $\ell \Subset \Omega$ such that $\P[\sle(\sigma)\in \ell]=0$. 
	Let $\lambda=\inf\{t\ge\sigma: \sle(t) \in \ell \}$ and $\lambda_\mesh=\inf\{t\ge \sigma_\mesh: \sle_\mesh(t)\in \ell  \}$. Then   $\lambda_\mesh\to \lambda$ a.s.  
\end{lemma}
\begin{proof}
	With probability $1$, there exist sequences of rational times $t_\eta\downarrow \lambda$ and $s_\eta\uparrow \lambda$ for $\eta\to 0$ in a countable set such that $\sle([s_\eta,t_\eta])\cap \ell\neq \emptyset$. This can be easily proved by way of contradiction, by using that an SLE$_6$ curve will a.s.\ cross a deterministic smooth curve upon hitting it.  Now the lemma follows from the continuity of $\gamma$.
\end{proof}

\subsection{Arm events}\label{subsec:arm}
Given a percolation configuration $\omega_\mesh$ and  an annulus $A$, we say that an \notion{alternating  4-arm  event} occurs for $A$ if and only if there are four disjoint monochromatic  paths connecting $\p_1 A$ and $\p_2 A$ such that the color sequence of the four paths is alternating between black and white.  There is an ambiguity in the definition due to the lattice effect at the boundary. However the precise convention does not matter as $\mesh\to 0$ so we ignore it. 
In the continuum, suppose $A$ is an annulus such that $\p_1 A$ and $\p_2 A$ are piecewise smooth.  For $A\Subset \Omega$,  a quad-crossing configuration $\omega\in\cH_\Omega$ is said to belong to the alternating   4-arm event of $A$  if there exist  quads $\Quad_i \subset \cQ_\Omega$, $i=1,2,3,4$, with the following properties:
\begin{itemize}
	\item[(i)] $Q_1$ and $Q_3$ are disjoint and at positive distance  from each other, and the same hold for $Q_2$ and $Q_4$.
	\item[(ii)] For $i\in\{1,3\}$, the side $\p_1Q_i$ lies inside $\p_1 A$ and the side $\p_3Q_i$ lies outside  $\p_2 A$; for $i\in\{2,4\}$, the side $\p_2Q_i$ lies inside $\p_1 A$ and the side $\p_4Q_i$ lies outside  $\p_2 A$; all these sides are of positive distance  away from $A$ and from the other $Q_j$'s.
	\item[(iii)] The four quads are ordered cyclically around $A$ according to their indices.
	\item[(iv)] $\omega(Q_1)=\omega(Q_3)=1$ and $\omega(Q_2)=\omega(Q_4)=0$.
\end{itemize}
In both the discrete and the continuum, the general $k$-arm event in $A$ given any prescribed color pattern can be defined similarly. For $\omega_\mesh$ coming from a percolation configuration, the two  definitions of arm events agree.

\begin{conv}\label{conv:arm}
In the rest of the  paper,  for each $k=2,3,4,5$, we focus on arm events with particular color conditions.  For $k=4$, it is the alternating 4-arm event. For $k=2,3,5$, it is the   $k$-arm event where not all arms have the same color.  
We will call  these events  the $k$-arm event  without mentioning the color pattern. We will not need the case $k\neq 2,3,4,5$.
\end{conv} 
Now we are ready to describe the normalizing constants in  \eqref{eq:im} and \eqref{eq:pm}.
\begin{remark}[Normalizing constants]\label{rmk:normalization}
	We use $\alpha^\mesh_k(\mesh,1)$ ($k=2,4$) 
	to denote the probability of the $k$-arm event (under Convention~\ref{conv:arm}) 
	from the single site at the origin to $\partial \bx_1$. Then $\alpha^\mesh_2(\mesh,1)$ and $\alpha^\mesh_4(\mesh,1)$ are the normalizing constants in \eqref{eq:im} and \eqref{eq:pm}, respectively. It is known that \(\alpha^\eta_k(\eta,1)=\mesh^{(k^2-1)/12+o_\mesh(1)}\) \cite{Sminov-Werner}. The up-to-constant asymptotics are open.
\end{remark}
In the coupling of Theorem~\ref{thm:embedding}, for $k=2,3,4,5$, let $\cA_k$ be $k$-arm events for an annulus $A\subset\Omega$ as in Convention~\ref{conv:arm}. Then the event
 $\cA_k$ is a.s.\ measurable with respect to the Borel $\sigma$-algebra of $(\cH_\Omega,d_\cH)$ \cite[Section~2]{GPS-piv}. 
As explained in \cite{Sminov-Werner},  the events  $\cA_k$ can be expressed in terms of percolation exploration to give
\begin{equation}\label{eq:arm-lim}
\lim_{\mesh\to 0} \P[\omega_\mesh\in \cA_k]=\P[\omega\in \cA_k] \qquad \textrm{for }k=2,3,4,5.
\end{equation}
Lemma 2.9  in \cite{GPS-piv} gives the following stronger version of \eqref{eq:arm-lim} when $k=2,3,4$.  (This is expected to be true also for $k=5$, but this is not proved in \cite{GPS-piv} and is not needed.)
\begin{lemma}\label{lem:arm}
$\lim_{\mesh\to0}\P[\{\omega\in \cA_k\}\xor \{\omega_\mesh\in \cA_k\}]=0$ for $k=2,3,4$.
\end{lemma}
For $R>r>0$ and $A=A(r,R)$ write $\alpha_k(r,R)=\P[\omega\in \cA_k]$. 
The up-to-constant asymptotic for $\alpha_k(r,R)$  is well-known \cite[Equation~(14)]{Sminov-Werner} (c.f.\ Remark~\ref{rmk:normalization}):
\begin{align}\label{eq:24'}
\alpha_k(r,R)\asymp(r/R)^{(k^2-1)/12}, \quad \textrm{for } k=2,3,4,5.
\end{align}

An important property of  $\omega$ as an element in $\cH$ is the monotonicity built in its definition. The following monotonicity results will be used repeatedly.
\begin{lemma}\label{lem:mono}
In the coupling of Theorem~\ref{thm:embedding}, let $\sle\defeq\sle^{ab}$ and $\sle_\mesh\defeq\sle^{ab}_\mesh$ for two given distinct points $a,b\in \partial \Omega$. View $(\sle_\mesh,\sle)$ as parametrized curves as in Lemma~\ref{lem:stop}.  For each fixed $t\in(0,1)$, let $K_t$ be the hull of $\sle([0,t])$. Namely, $K_t$ is the complement of the connected component of $\Omega\setminus \sle([0,t])$ containing the target of $\sle$. 
For any annulus $A\Subset \Omega$, let $A_1$ be the inside of $\p_1 A$ and $A_2$ be the outside of $\p_2 A $. 
Then for the quad-crossing configuration $\omega$,
\begin{enumerate}
			\item \label{item:2-arm}
			if $\p_1 A\cap \sle \neq \emptyset$, then the 2-arm event for $A$ occurs a.s., and
			\item \label{item:3-arm}
			if there exists $t\in[0,1]$ such that $\p K_t  \cap A_1 \neq\emptyset$ and  $\sle(t) \in A_2$, then   
			the 3-arm event for $A$ occurs a.s. 
\end{enumerate}	
\end{lemma} 
\begin{figure}
		\centering
		\includegraphics[scale=1]{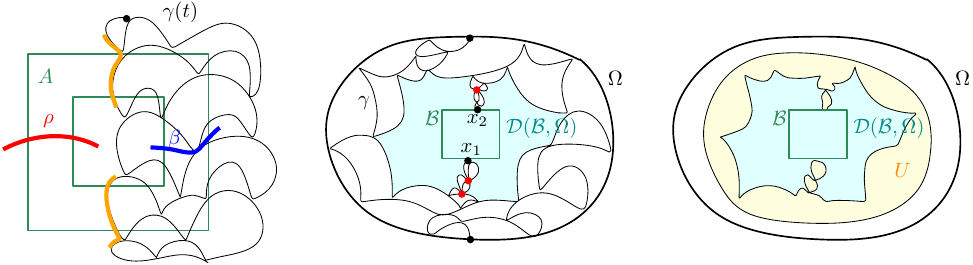}
	\caption{Left: Illustration of Lemma \ref{lem:mono}, Assertion 2. The arms are shown in blue and orange instead of black and white, respectively. Middle: Illustration in light blue of the face $\cD(\cB,\Omega)$ at $\mathcal B$ induced by $\gamma$ as defined in Section \ref{subsec:couple}. Double points of $\partial \cD(\cB,\Omega)$ are marked in red and correspond to local cut points for $\gamma$ (see Remark \ref{rmk:boundary}). {Right: Illustration of the event $\cG_\Omega(\bx,U)$.} } 
	\label{fig1}
\end{figure}
	\begin{proof}
		For Assertion 1, choose $t,\delta \in(0,1)$ such that $\sle(t)$ is inside $\p_1 A$ and $\dist(\sle(t), \p_1 A)>\delta$. 
		Since $\lim_{\mesh\to 0} \sup_{0\le t\le 1}\{ |\sle_\mesh(t) -\sle(t)|\}=0$ a.s., we  have  $\dist(\sle_\mesh(t), \p_1 A)>0.5\delta$ for small enough $\mesh$. In this case the 2-arm event for $A$ occurs for $\omega_\mesh$. Sending $\eta\to 0$ and applying Lemma~\ref{lem:arm}, we get Assertion~\ref{item:2-arm}.
		
		Assertion~\ref{item:3-arm} can be proved similarly. See Figure \ref{fig1}. Since $\sle(t)$ is a boundary point of $K_t$, there exists a $\delta>0$  and a  path $\rho$ starting from some point in  $A_1$ and ending at some point in  $A_2$ such that $ \dist(\rho, \sle([0,t])) > \delta$. Without loss of generality assume the set $\partial K_t\cap A_1$ contains a point on left frontier of $\gamma([0,t])$. Now for small enough $\mesh$,  we have $\dist(\rho, \sle_\mesh([0,t])) > 0.5\delta$. On the other hand, by the argument for Assertion~\ref{item:2-arm}, for $\mesh$ small enough the there exists a black arm $\beta$ of $\omega_\mesh$ from $\p_1 A$ to  $\p_2 A$. In this case, there must be a white arm of $\omega_\mesh$ on each connected component of $A\setminus (\rho \cup \beta)$ from $\p_1 A$ to  $\p_2 A$, hence the 3-arm event for $A$ occurs for $\omega_\mesh$. Now Assertion~\ref{item:3-arm} follows from Lemma~\ref{lem:arm}.
	\end{proof}
The following variant of Lemma~\ref{lem:mono} can be proved similarly. We omit the details.
	\begin{lemma}\label{lem:mono2}
Consider a coupling where both the conditions in Theorems~\ref{thm:length2} and~\ref{thm:embedding} are satisfied so that $(\omega,\sle^1,\sle^2)$ are coupled. 
Let $A,A_1,A_2$ be defined as in Lemma~\ref{lem:mono}.  Then  on the event $\cP\cap A_1\neq \emptyset$,  the 4-arm event for $A$ occurs a.s. for $\omega$.
	\end{lemma}
The event $\gamma\cap A_1\neq\emptyset$ in Lemma~\ref{lem:mono} is simply the 2-arm event with the further requirement that each of the two boundary arcs contain one endpoint of the arm.
The similar statement holds for $\cP\cap A_1\neq \emptyset$ in Lemma~\ref{lem:mono2}. 
By the following lemma, these endpoint requirements only decrease the probability by a constant factor.
\begin{lemma}\label{lem:endpoint}
In the setting of Lemmas~\ref{lem:mono} and~\ref{lem:mono2}, let $\bx\subset \Omega$ 
be a box of radius $\eps$  whose center is $r>10\eps$ away from $\p\Omega$. Then 
$\P[\sle\cap \bx\neq \emptyset]\asymp \alpha_2(\eps,r)$ and $\P[\cP\cap \bx\neq \emptyset]\asymp \alpha_4(\eps,r)$, where the implicit constants in $\asymp$ only depend on $\Omega$ but  not on other parameters.
\end{lemma}
\begin{proof}
The 2-arm case follows from the classical one-point estimate of $\SLE_6$. See e.g.\ \cite{Lawler-Minkowski}. 
The 4-arm case  follows from \cite[Proposition~4.9]{GPS-piv}.
\end{proof}

\subsection{Face induced by the percolation exploration}\label{subsec:couple}
Given a box $\bx$ and two distinct points $x_1,x_2\in \p \bx$,
let $\theta_1$ (resp., $\theta_2$) be a simple path joining $x_1$ and $x_{2}$ (resp., $x_2$ and $x_1$).  
	If the pair of  paths
	$\Theta=\{\theta_1,\theta_2 \}$ is such that there exists a domain $\cD_\Theta$  with $\bx\subset \cD_\Theta$ and  $\p \cD_\Theta=\theta_1\cup\theta_2$,
	then we call	$\Theta$  a   \notion{face} at $\bx$
	with endpoints $x_1,x_2$. 
	
	Let $\Omega$ be a simply connected domain whose boundary is a continuous curve and let $a, b\in \p \Omega$ be such that $a \neq b$.  
	Suppose $\sle$ is a SLE$_6$ on $(\Omega,a,b)$ parametrized in an arbitrary way and $\bx \subset \Omega$ is a box.
	Throughout this subsection we write $(\Omega,a,b)$ as $\Omega$ whenever it simplifies the notation and cause no confusion. For example, we 
	use $\cA(\bx,\Omega)$ to denote the event $\{\sle\cap\bx \neq \emptyset\}$ although this event depends on $a,b$. On $\cA(\bx,\Omega)$, let 
	\begin{align*}
		&\ubar{\sigma}=\inf\{t:\sle_t\in \bx\}, \quad&& \ol\sigma=\sup\{t:\sle_t\in \bx\},\\
		&x_1=\sle(\ubar{\sigma}), &&x_2=\sle(\ol\sigma).
	\end{align*}
	Let $\mathcal{D}(\bx,\Omega)$ be the connected component of $\Omega  \setminus (\sle[0,\ubar\sigma]\cup   \sle[\ol \sigma,\infty) )$ containing $\bx$. 
	Then $\cD(\bx,\Omega)$ can be viewed as a (random) face  at $\bx$ with  the arcs $\p_{x_1,x_2}\cD(\bx,\Omega)$ and $\p_{x_2,x_1}\cD(\bx,\Omega)$,
	which we call \notion{the face at  $\bx$ induced by $\sle$}.	 By setting $\cD(\bx,\Omega)=\emptyset$ when $\cA(\bx,\Omega)$ does not occur,  we can view $\cD(\bx,\Omega)$  as a random domain with  two ordered boundary marked points $x_1,x_2$ when it is nonempty. Given a simply connected domain $U$  with piecewise smooth boundary such that  $\bx\Subset U\Subset \Omega$, set
	\begin{equation}\label{eq:Gevent}
		\cG_\Omega(\bx,U)\defeq\cA(\bx,\Omega)\cap\{\cD(\bx,\Omega)\subset U \}.
	\end{equation}
	{See the right part of Figure \ref{fig1} for an illustration.} The picture above has a discrete counterpart. Suppose $\Omega$ is a Jordan domain and $\omega_\mesh$ is sampled from $\Ber(\Omega_\mesh)$. Let $\sle_\mesh$ denote the associated  interface on $(\Omega_\mesh,a_\mesh,b_\mesh)$ for some $a\neq b\in \p \Omega$. 
	Let $\cA_\mesh(\bx,\Omega)$ be the event that there exists an edge on $\gamma_\mesh$ such that the two hexagons containing the edge are both in $\bx$. 
	Consider the first and last such edges on $\sle_\mesh$, whose visiting time are denoted by $\ul\sigma_\mesh$ and $\ol\sigma_\mesh$, respectively. Let $\Theta_\mesh$ denote the face at $\bx$ induced by $\sle_\mesh$, which forms the boundary of the domain $\cD_\mesh(\bx,\Omega)$. Similarly as in \eqref{eq:Gevent}, define $\cG_{\Omega_\eta}(\bx,U)\defeq\cA_\mesh(\bx,\Omega) \cap  \{\cD_\mesh(\bx,\Omega)\subset U\}$.  We have the following two lemmas.
	\begin{lemma}\label{lem:GG}
		Suppose $\Omega$ is a Jordan domain and that $\bx$ and $U$ are defined as above. Suppose we are in the coupling of Theorem~\ref{thm:embedding}. We view
		 $\sle_\mesh$ and $\sle$ as parametrized curves as in Lemma~\ref{lem:stop}.  Then
		$\ol{\cD_\mesh(\bx,\Omega)}$ converges to $\ol{\cD(\bx,\Omega)}$ in probability for the Hausdorff metric as closed sets with two ordered marked points. 
	\end{lemma}
\begin{proof}
It suffices to show that
\begin{equation}\label{eq:box}
\lim_{\mesh \to 0}\P[ \{\cB' \cap \ol{\cD_\mesh(\bx,\Omega)} =\emptyset \} \Delta \{\cB' \cap  \ol{\cD(\bx,\Omega)}=\emptyset \}] =0,\quad \textrm{for a fixed box}\; \bx'\Subset \Omega.
\end{equation}	
Given a fixed piecewise smooth curve $p:[0,1]\to \Omega$ with $p(0)\in \p \Omega, p(1) \in \bx$ and $p((0,1))\subset \Omega$. If $p\cap \p \cD(\bx,\Omega)=\emptyset$, since $\sle_\mesh$ converges to $\sle$ in the  $d_\cU$-metric, for small enough $\mesh$ we must have $p\cap \cD_\mesh(\bx,\Omega)=\emptyset$. If $p \cap \p \cD(\bx,\Omega)\neq \emptyset$, then by Lemma~\ref{lem:stop}  for small enough $\mesh$ we must have $p\cap\cD_\mesh(\bx,\Omega)\neq \emptyset$.  This implies \eqref{eq:box} by elementary topological consideration.
\end{proof}
\begin{lemma}\label{lem:GG2}
In the setting of Lemma~\ref{lem:GG}, let $\cA_3$ represent the $3$-arm event for $U\setminus \bx$. Then    \(\Pb{\cA (\bx, \Omega)\setminus\cG_{\Omega}(\bx,U)} \le \P[{\omega\in \cA_3}].\)
\end{lemma}
\begin{proof}	
By Lemma~\ref{lem:stop}, $\P[\ol{\cD(\bx,\Omega)}\subset U]=\P[\ol{\cD(\bx,\Omega)}\subset \ol U]$.
 By  Lemma~\ref{lem:GG} it suffices to show that \(\cA_\mesh(\bx, \Omega)\cap\{\omega_\mesh\notin \cA_3 \} \subset \cG_{\Omega_\mesh}(\bx,U).\)  To prove this, we see that  if $\cA_\mesh (\bx, \Omega)\cap \{\omega_\mesh\notin \cA_3\}$ occurs,  the black sites adjacent to  $\sle_\mesh ([0,\underline\sigma_\mesh])$ and $\sle_\mesh([0,\ol \sigma_\mesh])$ must share a common hexagon within $U\setminus \bx$. The similar statement holds for the white sites. This concludes the proof.
\end{proof}
In the setting of Theorem~\ref{thm:embedding},  it is clear that $\omega$ inherits the  spatial independence property from $\omega_\mesh$. By Lemma~\ref{lem:GG}, we get the following.
	\begin{lemma}\label{lem:independent}
In the setting of Lemma~\ref{lem:GG}, let  $\Omega_1,\Omega_2$ be two disjoint open subsets in $\Omega$.
Then $\omega$ restricted to $\cQ_{\Omega_1}$ and to $\cQ_{\Omega_2}$ are independent as random variables in $\cH_{\Omega_1}$ and $\cH_{\Omega_2}$, respectively.  Moreover,  $\omega$ restricted to   $\cQ_\cB$ is independent of  $\cD(\bx,\Omega)$.
\end{lemma}

\section{Equivalence of the two measures on the interface}\label{sec:equiv}
This section is devoted to proving Theorem~\ref{thm:equi-length} hence we retain the notations in the statement of the theorem. 
{To prove Theorem~\ref{thm:equi-length},  we use the $L^2$ framework as in \cite{GPS-piv}, which is based on a strong coupling scheme and the spatial independence  of percolation. 
Since we work in the continuum, some issues in \cite[Section 4]{GPS-piv} can be simplified. In particular, the required one-point and two-point estimates  that we will rely on are 
are power laws with no sub-polynomial corrections (see Lemmas~\ref{lem:12y}---\ref{lem:two-arm}), 
while a major novelty of  \cite[Section 4]{GPS-piv} is obtaining scaling limit results despite the unknown sub-polynomial corrections in the percolation estimates. 
After we prepare the one-point and two-point estimates, we reduce Theorem~\ref{thm:equi-length} to a strong coupling estimate \eqref{eq:coup}. This reduction is a straightforward adaptation of the $L^2$ argument in \cite[Section 4]{GPS-piv}, nevertheless we still include the full argument for completeness and hence follow closely both the  method and the presentation in \cite[Section 4]{GPS-piv}.
To prove the strong coupling estimate~\eqref{eq:coup}, we would like to apply its discrete analog from \cite{GPS-piv} and then pass to the continuum. However, a straightforward implementation of this idea only gives Lemma~\ref{prop:coupling}, a weaker variant of~\eqref{eq:coup}. The reason is that when we pass from percolation to its continuum limit, we rely on Theorem~\ref{thm:embedding}, which is for Jordan domains. On the other hand, the domain boundary considered in~\eqref{eq:coup} is the exterior boundary of $\SLE_6$, which is not simple. Instead of trying to strengthen the convergence in Theorem~\ref{thm:embedding} to include certain non-Jordan domains, we will use an argument directly in the continuum to go from Lemma~\ref{prop:coupling} to the desired~\eqref{eq:coup}. We now carry out the plan above in detail. 
}

Let $\cB\Subset\Omega$ be a box whose four vertices are on $\cup_{k\in \N} 2^{-k}\Z^2$. Let $\eps\in \{2^{-k}: k\in \N \}$.
Assume $\eps$ is small enough such that $\eps<\dist(\bx, \p \Omega)$. 
Then $\bx$ is partitioned by certain boxes of radius $\eps$ centered at points on the lattice $2\eps \,\Z^2$.
Let $Q_1, \dots, Q_p$ be a list of these boxes in arbitrary order. For $i\ge 1$, let $q_i$ denote the center of $Q_i$. 	Let 
$$
Y^\eps=\# \{1\le i \le p:   \sle\cap 2Q_i\neq \emptyset \}.
$$
In \cite[Section 5.3]{GPS-piv}, the following is proved.%
\footnote{To obtain Proposition \ref{prop:GPS2} from \cite[Section 5.3]{GPS-piv} we use that, in the notation of that paper, $X$ appropriately renormalized converges to $\tau$, $\E[ (X-\beta_{\mathrm{two-arm}} Y)^2 ]=o(\E[X^2])$, and $\beta_{\mathrm{two-arm}}\asymp \eps^2\eta^{-2}\alpha_4^\eta(\eta,\eps)$.}
\begin{proposition}\label{prop:GPS2}
	There exists a deterministic constant $c>0$ such that
	\[
	\tau(\bx) =\lim_{\eps\to 0} \frac{c \, Y^{\eps}}{\eps^{-2} \alpha_2(\eps ,1)}
	\qquad\hbox{in $L^2$}.
	\]
\end{proposition}
Consider  the square $Q_0\defeq \bx_\eps(0)$.
Let $\sle^0$ be a chordal SLE$_6$ on  $(\cB_1, -i,i)$  and $x_0\defeq\omo(Q_0)$, where $\omo$ is the occupation measure of $\sle^0$. 
Let $\mathcal{A}_0(2\eps,1)$ be the event that $\sle^0\cap 2Q_0\neq \emptyset$, and define
\begin{equation}\label{e.beta}
	\nor\defeq\Eb{ x_0 \bigm| \mathcal{A}_0(2\eps,1)}.
\end{equation}	
\noindent Theorem~\ref{thm:equi-length} is an immediate consequence of Proposition~\ref{prop:GPS2} and the following.
\begin{proposition}\label{l.XY}
For each {box} $\cB\Subset\Omega$ {as above}, we have that $\om(\cB) =\lim_{\eps\to 0}\nor Y^\eps $ in $L^2$.
\end{proposition}	

{Before proving Proposition~\ref{l.XY}, we first record a few basic estimates in Lemmas~\ref{lem:12y}---\ref{lem:two-arm}.}

Define $y_i$ to be the indicator function of the event that $\sle\cap 2Q_i\neq \emptyset$ so that 
$Y^\eps= \sum_{1}^{p}y_i$. Similarly, for any $1\le i \le p$, let $x_i=\om(Q_i)$ such that  $	\om (\bx)=\sum_{1}^{p}x_i$.
We first record some a priori estimates for the $x_i'$s and the $y_i$'s. These estimates would trivially follow from known Green function estimates for $\SLE_6$ \cite{Lawler-Green}. 
However,  we instead present an argument   that can be  readily extended to the case of pivotal points in Section~\ref{sec:piv}. {The following result is classical, and we refer to \cite{beffara-dim} for a proof.}
\begin{lemma}\label{lem:12y}
	In the above setting, for all $1\le i,j \le p$ with $i\neq j$, 
	\begin{equation}\label{eq:1pty}
		\E[y_i] \asymp \eps^{1/4} \qquad\textrm{and}\qquad	\Eb{y_i y_j}\lesssim \frac{\eps^{1/2} }{|q_i-q_j|^{1/4}}.
	\end{equation}
	where the constants in $ \asymp$ and  $\lesssim$ only depend on $\cB$ and $\Omega$.
\end{lemma} 

A similar argument based on arm exponents gives the following.
\begin{lemma}\label{lem:12pt}
	For all $1\le i,j \le p$ with $i\neq j$, we have
	\begin{equation}\label{eq:2pt}
		\E[x_i]\lesssim \eps^2,\qquad \Eb{x_i x_j}\lesssim \frac{\eps^4 }{|q_i-q_j|^{1/4}} \qquad \textrm{and} \qquad\E[x^2_i]\lesssim \eps^{15/4},
	\end{equation}
	where  the constants in $\lesssim$ only depend on $B$ and $\Omega$.
\end{lemma}
\begin{proof}
	For $r\in (0,0.01\eps)$ and $\bullet=i,j$, let $\cX_\bullet=\sle\cap 2Q_\bullet$ and $\cX_\bullet^r=\{z: \dist(z,\cX_\bullet)\le r \}$.
	It is clear that $\cX^r_i\subset 4Q_i$. By Lemma~\ref{lem:mono} and \eqref{eq:24'}, $\P [\dist(z, \cX_i)\le r]\lesssim r^{1/4}$ for all $z\in 4Q_i$. 
	Therefore, by Fubini's theorem, we have 
	\[
	\E [\Area(\cX^r_i)] =\int_{4Q_i} \P [z\in \cX^r_i] dz \lesssim \eps^2 r^{1/4}.
	\]
	Now  Fatou's lemma and Definitions~\ref{def:Mink} and \ref{def:occu} yield $\E[x_i]\lesssim \eps^2$.
	
	For the second inequality, by Fubini's theorem, we have 
	\[
	\E[\Area(\cX^r_i)\Area(\cX^r_j)] =\int_{4Q_i\times 4Q_j} \P [z\in \cX^r_i, w\in \cX^r_j] \,dz\,dw.
	\]
	{By  Lemma~\ref{lem:12y},} we have $ \P [z\in \cX^r_i, w\in \cX^r_j] \lesssim r^{1/2}/|z-w|^{1/4}$. 
	Now the second inequality follows from Fatou's lemma and Definitions~\ref{def:Mink} and \ref{def:occu}.
	
	The third inequality follows from a similar argument as for the second one.
\end{proof}

By Lemmas~\ref{lem:12y} and~\ref{lem:12pt}, we have
\begin{equation}\label{eq:beta}
	\nor \le \frac{\Eb{x_0}}{\P[\mathcal{A}_0(2\eps,1)]} \lesssim \frac{\eps^2}{\eps^{1/4}}= \eps^{7/4}.
\end{equation}	
\begin{lemma}\label{lem:two-arm}
	In the above setting, for all $1\le i\le p$, let 
	\[
	\wt{\cX}_i =\cap_{\delta>0} \{z\in Q_i : \textrm{the 2-arm event occurs for the annulus $\bx(q_i, \tfrac32\eps)\setminus \bx(z,\delta)$} \}
	\]
	and $\wt \cX_i^r=\{z: \dist(z,\wt \cX_i)\le r \}$. Let  $\wt{x}_i = \liminf_{r \to 0} r^{-1/4}\Area (\wt\cX_i^r)$. Then
	\[
	x_i\le {\wt{x}_i} \quad \textrm{and}\quad \Eb{{\wt{x}_i}}\lesssim \eps^{7/4}\qquad\textrm{for all}\qquad 1\le i\le p
	\]	where the constant in $\lesssim$ is independent $\eps,i,\bx,\Omega$.
\end{lemma}	 
\begin{proof}
	Lemma~\ref{lem:mono} and \eqref{eq:Mink} imply that $x_i\le {\wt{x}_i}$. The bound $ \Eb{{\wt{x}_i}}\lesssim \eps^{7/4}$ follows from the same argument as for the first inequality in Lemma~\ref{lem:12pt}. Here the domain $\Omega$ is replaced by $\tfrac32 Q_i$. Therefore we get the upper bound $\eps^{7/4}$ instead of $\eps ^2$.
\end{proof}	
\noindent The advantage of considering ${\wt{x}_i}$ instead of $x_i$ is that it is completely determined by $\omega$ restricted to $\bx(q_i, \tfrac32\eps)$,  hence is independent of what happen outside   $2Q_i$.

\medskip

Now we proceed to prove  Proposition~\ref{l.XY}. Fix some $r>0$ to  be determined later. Write $\Delta_i=x_i - \nor y_i$ for $1\le i\le p $ and 
\begin{equation*}
	\Eb{(\om(\bx) - \nor Y^\eps)^2}= \sum_{i,j=1}^p \Eb{\Delta_i\Delta_j}.
\end{equation*}
Split the summation into an ``on-diagonal''  term and an ``off-diagonal'' term:
\begin{align}
	\label{e.correls2}
	\Eb{(\om(\bx) - \nor Y^\eps)^2} = \sum_{|q_i-q_j|\le r}\Eb{\Delta_i\Delta_j}
	+ \sum_{|q_i-q_j|>r}
\Eb{\Delta_i\Delta_j}.
\end{align}
To estimate the on-diagonal term, take any $i,j$ such that $|q_i-q_j|\le r$, and observe that since all variables and constants are positive, we have
\begin{equation}
	\label{e.xixj}
\Eb{\Delta_i\Delta_j}\le \Eb{x_i x_j +\nor^2 y_i y_j}.
\end{equation}
There are $O(1) \eps^{-2}$ choices for the box $Q_i$ (where $O(1)$ depends on $\cB$). For a fixed box $Q_i$ and any $k\geq 0$ such that $2^k \eps < r $, there are $O(1) 2^{2k}$ boxes
$Q_j$ satisfying $2^k \eps\le |q_i-q_j|<2^{k+1} \eps$. For any of these boxes, Lemma~\ref{lem:12pt} gives \(\Eb{x_i x_j}\lesssim \eps^4 /(2^k \eps)^{1/4}\).

Therefore
\begin{equation}\label{eq:xixj}
	\sum_{|q_i-q_j|\le r} \Eb{x_i x_j} \lesssim \eps^{-2} \sum_{k\le \log_2(r/\eps)} 2^{2k} \cdot
	\frac{\eps^4 }{(2^k \eps)^{1/4}}.
\end{equation}
By Lemma~\ref{lem:12y} and \eqref{eq:beta} we obtain the same bound on $	\sum_{|q_i-q_j|\le r} \Eb{\nor^2 y_iy_j}$.
Therefore 
\begin{equation}\label{eq:dia}
	\sum_{|q_i-q_j|\le r}
\Eb{\Delta_i\Delta_j} \lesssim r^{7/4}.
\end{equation}	

Now consider the off-diagonal term in \eqref{e.correls2}. We claim that for fixed $\delta$, if $\eps$ is small enough, for any $i,j$ such that  $l\defeq|q_i-q_j|>r$ we have
\begin{equation}\label{eq:two-pt}
	\Eb{\Delta_i\Delta_j} \le \delta \cdot  \frac{\eps^4}{l^{1/4}}.
\end{equation}

Let $\intd  \in (2\eps,r/4)$ be some intermediate distance whose value will be fixed later. For $k=2,3$ and $\bullet=i,j$, let $\cA^\bullet_k
=\mathcal{A}^\bullet_k(\intd, l/2)$  be the $k$-arm event for the annulus
$\bx(q_\bullet, l/2) \setminus \bx(q_\bullet, \intd)$.
Following the notations of Section~\ref{subsec:couple}, let $\cD_\bullet\defeq\cD(\bx(q_\bullet,\intd),\Omega)$ and $\Theta_\bullet$ be the face at $ \bx(q_\bullet,\intd)$ induced by 
$\gamma$. Let $\cG_\bullet=\cG(\bx(q_\bullet,\intd), \bx(q_\bullet, l/2))$. Note that by Lemma~\ref{lem:mono}, we have $\cG_\bullet\subset \cA^\bullet_2$.

Let $\cW=\cG_i\cap \cG_j$ and $\cZ=(\cA^i_2\cap \cA^j_2) \setminus \cW$.		By Lemma~\ref{lem:mono}, if $\Delta_i\Delta_j\neq 0$, the event $\A^i_2\cap \cA^j_2$ must occur.  Therefore 
\begin{align}
	\label{e.ZiWi}
	\Eb{\Delta_i\Delta_j} = \Eb{\Delta_i\Delta_j \1_{\cZ}} + \Eb{\Delta_i\Delta_j \1_{\cW}}.
\end{align}
Let $\mathcal{A}_{i,j}$ be the event that two-arm events occur in the annuli  $\bx(q_j, l/2)\setminus 2Q_j$,
$\bx(q_i, l/2)\setminus 2Q_i$   and $\Omega\setminus  \bx(\frac {q_i+q_j} 2, l)$. Observe that if $(x_i x_j + \nor^2 y_i y_j)\neq 0$ then $\mathcal{A}_{i,j}$ occurs. Recall ${\wt{x}_i}$ in Lemma~\ref{lem:two-arm}.  We have
\[
\Eb{|\Delta_i\Delta_j| \1_{\cA^i_2\setminus \cG_i}} 
\le 
\Eb{(x_i x_j + \nor^2 y_i y_j) \1_{\cA^i_2\setminus \cG_i} }
\le 	\Eb{({\wt{x}_i}{\wt{x}_j} +\nor^2  )\cdot  \1_{\cA^i_2\setminus \cG_i} \cdot \1_{\mathcal{A}_{i,j}} }.
\]	
By Lemma~\ref{lem:independent}, ${\wt{x}_i}$, ${\wt{x}_j}$ and $ \1_{\cA^i_2\setminus \cG_i} \cdot \1_{\mathcal{A}_{i,j}} $ are  independent.
By Lemma~\ref{lem:two-arm}  and \eqref{eq:beta}, we have
\[
\Eb{({\wt{x}_i}{\wt{x}_j} +\nor^2  )\cdot  \1_{\cA^i_2\setminus \cG_i} \cdot \1_{\mathcal{A}_{i,j}} }=\left(\Eb{{\wt{x}_i}}\Eb{{\wt{x}_j}} +\nor^2\right)\P[(\cA^i_2\setminus \cG_i)\cap \mathcal{A}_{i,j}]
\lesssim \eps^{7/2} 	\P[(\cA^i_2\setminus \cG_i)\cap \mathcal{A}_{i,j}].
\]
By the same argument as in Lemma~\ref{lem:GG2}, we have \(\P[(\cA^i_2\setminus \cG_i)\cap \mathcal{A}_{i,j}] \le\P[\omega\in\cA^i_3\cap \mathcal{A}_{i,j}]\). 
By Lemma~\ref{lem:independent} and~\eqref{eq:24'}, we have $\P[\omega\in\cA^i_3\cap \mathcal{A}_{i,j}]=o_{\intd/l}(1) \eps^{1/2}/l^{1/4}$. Therefore 
\[
\Eb{|\Delta_i\Delta_j| \1_{\cA^i_2\setminus \cG_i}}=o_{\intd/l}(1) \frac{\eps^4}{l^{1/4}}.
\]
Since $\cZ\subset (\cA^i_2\setminus \cG_i) \cup (\cA^j_2\setminus \cG_j)$, we have 
\begin{equation}\label{eq:bad}
	\left|	\Eb{\Delta_i\Delta_j  \1_{\cZ}} \right| =o_{\intd/l}(1) \frac{\eps^4}{l^{1/4}}.
\end{equation}
\vskip 0.3 cm
It remains to bound the second term on the right side of~\eqref{e.ZiWi}.  Recall the notations introduced in Section~\ref{subsec:couple}.  
For $\bullet=i,j$, on the event $\cG_\bullet$,	let $\ubar{\sigma}_\bullet=\inf\{t:\sle_t\in B(q_\bullet,\intd)\}$ and $\ol\sigma_\bullet=\sup\{t:\sle_t\in B(q_\bullet,\intd)\}$.  
Let $\sle^\bullet$ be the curve $\sle([\ubar{\sigma}_\bullet,\ol\sigma_\bullet])$. Then by the reversibility of SLE$_6$, the curve $\sle^\bullet$  conditioning on $\cD_\bullet$ is a chordal SLE$_6$ inside $\cD_\bullet$.
We claim that
\begin{equation}\label{eq:coup}
	\1_\cW \big|\Eb{x_i - \nor y_i \bigm| \cD_{i}}\big| =o_{\eps/\intd} (1) \frac{\eps^2}{\intd^{1/4}} \quad \textrm{and the same with $j$ in place of $i$}.
\end{equation}
Let us first wrap up the proof of Proposition~\ref{l.XY} given \eqref{eq:coup}. On $\cW$, the curves $\sle^i,\sle^j$ are independent conditioned on $\cD_i,\cD_j$.  Combining with \eqref{eq:coup}, we get
\[
\Eb{ \1_{\cW}\left| \Delta_i\Delta_j \right|}= o_{\eps/\intd} (1) \frac{\eps^4}{\intd^{1/2}} \P[\cW].
\]	
On $\cW$, the 2-arm event occurs in the disjoint  annuli $\Omega\setminus \bx(\frac{q_i+q_j} 2, l),  \bx(q_i, l/2) \setminus \bx(q_i,\intd)$ and $\bx(q_j, l/2) \setminus \bx(q_j,\intd)$.
By Lemma~\ref{lem:independent}, we have $\P[\cW] \lesssim \intd^{1/2}/l^{1/4}$.  Therefore,
\begin{eqnarray}
	\big|\Eb{\1_{\cW}\Delta_i\Delta_j }\big| = o_{\eps/\intd} (1) \frac{\eps^4}{l^{1/4}}.
\end{eqnarray}
Combining with \eqref{eq:bad} and setting $\intd=r^2=\eps^{1/2}$, we get \eqref{eq:two-pt}. Summing over $i,j$, we see that the off-diagonal term in \eqref{e.correls2} is less than $\delta$ for sufficiently small $\eps$. Combining with \eqref{eq:dia}, this concludes the proof of  Proposition~\ref{l.XY}, and hence of Theorem~\ref{thm:equi-length}.

Now we focus on the proof of  \eqref{eq:coup}, which crucially relies on the following lemma.
\begin{lemma}\label{prop:coupling}
	Let $\Omega'$ be a Jordan domain containing 0. Let $d=\dist(0,\partial \Omega')$ and $d'=d\wedge 1$. 
	Let $a',b'\in \partial \Omega$  and $\gamma'$ be a chordal $\SLE_6$ on $(\Omega,a',b')$. Let 
	\[
	x'=\Mink_{7/4}(\sle'\cap \bx_{2\eps})\qquad \textrm{and}\qquad y'=\1_{\sle'\cap \bx_{2\eps}\neq\emptyset}.
	\]
	Then there exist absolute constants $c,C>0$ independent of $\Omega'$ such that for $0<\eps <d'/10$, 
	\begin{align}\label{eq:couple}
		\left|\Eb{  x' -   \nor  y'}\right|  \le C\left(\frac {2\eps}  {d'}\right)^{c}  \cdot  \eps^{7/4} \cdot \alpha_2 (2\eps,d').
	\end{align}		
\end{lemma}
\begin{proof}
	Recall $ x_0$ and $\mathcal{A}_0(2\eps,1)$ in  the definition of  $\nor$. Also recall the notations in Section~\ref{subsec:couple}. Set $\cA\defeq\cA(\bx_{2\eps},\Omega')=\{\sle'\cap \bx_{2\eps}\neq\emptyset\}$. 
	We have
	$$
	\Eb{  x' -   \nor  y'\bigm|\mathcal{A}}  = \Eb{  x' \bigm|
		\cA} - \Eb{ x_0 \bigm| \mathcal{A}_0(2\eps,1)}.
	$$
	Suppose $\omega'_\mesh$ is a site  percolation configuration on $ \Omega'_\mesh\setminus \bx_{1.9\eps}$. Then the discrete analog $\cA_\mesh(\bx_{2\eps},\Omega)$ of $\cA$ is an event measurable with respect to $\omega'_\mesh$.  Moreover, the face $\cD_\mesh(\bx_{2\eps},\Omega')$ induced by $\sle'$ at $\bx_{2\eps}$ is also measurable with respect to $\omega'_\mesh$.  
	Now assume the law of $\omega'_\mesh$ is the critical percolation conditioning on $\cA_\mesh$.  Let $\omega^0_\mesh$ be the random site percolation configuration defined in the same manner as $\omega'_\mesh$ with $(\bx_1,-i,i)$ in place of $(\Omega',a',b')$. 
	
	By \cite[Proposition 3.6]{GPS-piv}, there exist an absolute constant $c>0$ independent of $\Omega'$ and  a coupling $(\omega'_\mesh,\omega^0_\mesh)$ such that for $10\mesh<\eps<d'/10$, with probability at least $1-(2\eps/d')^c$, 
	we have $\cD_\mesh(\bx_{2\eps},\Omega')=\cD_\mesh (\bx_{2\eps},\bx_1)$. In fact, \cite[Proposition 3.6]{GPS-piv} is stated for the 4-arm event but as explained in \cite[Section 5.3]{GPS-piv}, the result holds for the 2-arm case here with little adaption. In this coupling, we extend $\omega'_\mesh$ and $\omega^0_\mesh$ to $\bx_{1.9\eps}$ by coloring each vertex black with probability $1/2$ and 
	white with probability $1/2$. Here we use the same randomness for $\omega'_\mesh$ and $\omega^0_\mesh$ on $\bx_{1.9\eps}$ while  different vertices are colored independently. 
	By Theorem~\ref{thm:embedding} and Lemma~\ref{lem:GG}, letting $\mesh\to 0$, we have a continuum coupling $(\ol \gamma',\ol \gamma^0,\omega',\omega^0)$ such that 
	\begin{itemize}
		\item $\ol \gamma'$ and $\ol \gamma^0$ are the scaling limits of the interfaces of $\omega'_\mesh$ and $\omega^0_\mesh$, respectively;
		\item $\omega'$ and $\omega^0$ are the scaling limits of $\omega'_\mesh$ and $\omega^0_\mesh$, respectively, as quad-crossing configurations in the $d_\cH$ metric;
		\item $\omega'$ has the law of $\omega$ as  in Theorem~\ref{thm:embedding} with $\Omega'$ in place of $\Omega$, 
		conditioning on $\cA$;
		\item the law $\omega^0_\mesh$  is the same as $\omega'_\mesh$ with $(\bx_1,-i,i)$ in place of  $(\Omega',a',b')$;
		\item with probability at least $1- (2\eps/ d')^{c}$,  we have $\cD(\bx_{2\eps},\Omega')=\cD(\bx_{2\eps},\bx_1)$;
		\item $\omega'=\omega^0$ inside $\bx_{1.9\eps}$, which is independent of $\cD(\bx_{2\eps},\Omega')$ and $\cD(\bx_{2\eps},\bx_1)$.
	\end{itemize}
	Let $F$ be the event that $\{\cD(\bx_{2\eps},\Omega')=\cD(\bx_{2\eps},\bx_1)\}$. Let $\ol x',\ol x_0$ be defined in the same way as $x',x_0$ with $(\sle',\sle^0)$ replaced by $(\ol{\sle}',\ol{\sle}^0)$. (Here the only difference between  $(\sle',\sle^0)$ and $(\ol{\sle}',\ol{\sle}^0)$ is that the former is unconditioned and the latter is conditioned.)
	Then $  \ol x' = \ol x_0$ on $F$ and  $\P[F]\ge 1-(2\eps/d')^c$.
	Therefore 
	\begin{align}
		\label{e.hatxi}
		\left|\Eb{  x' -  \nor  y' \bigm| \cA}\right|  \le  \left(\frac {2\eps}  {d'}\right)^c
		\big(\Eb{  \ol x'  \bigm|F^c}
		+ \Eb{ \ol x_0 \bigm|F^c} \big).
	\end{align}
	Let  $\wt{x}'$ be defined as  in Lemma~\ref{lem:two-arm} with $\ol{\sle}'$ in place of $\gamma$. Then $\ol x'\le \wt{x}'$. By the nature of the coupling,   $\wt{x}'$ is independent of $F$. Therefore, 	
	\[
	\Eb{  \ol x'  \bigm|F^c}\le \Eb{\wt{x}'}\lesssim \eps^{7/4}.
	\]
	Similarly, we have $\Eb{ \ol x_0 \bigm| F^c} \lesssim \eps^{7/4}$.
	Combining with \eqref{e.hatxi}, and using that $x' -   \nor  y'=0$ unless $\mathcal A$ occurs, we see that there exists a constant $C>0$
	\begin{align*}
		\left|\Eb{  x' -   \nor  y'}\right|  \le C\left(\frac {\eps}  {d'}\right)^{c}  \cdot  \eps^{7/4} \cdot \P[\cA].
	\end{align*}	
	Now Lemma~\ref{lem:mono} yields \eqref{eq:couple}.
	\qedhere
\end{proof}
\begin{remark}\label{rmk:boundary}
	We assume that $\Omega'$ is a Jordan domain in Lemma~\ref{prop:coupling} because our proof crucially relies on the coupling result of \cite{GPS-piv} in the discrete and the convergence result Theorem~\ref{thm:embedding}, which  is only established for Jordan domains \cite{Camia-Newman-CLE1}. Lemma~\ref{prop:coupling} is not directly applicable to $\cD_i,\cD_j$ in \eqref{eq:coup} since they are a.s.\ not Jordan (see Figure \ref{fig1}). To overcome this issue, we extend Lemma~\ref{prop:coupling} to Lemma~\ref{lem:coupling2} below.
\end{remark}
\begin{lemma}\label{lem:coupling2}
	Suppose $\Omega'$ is a simply connected domain containing the origin whose boundary is a continuous curve.  Let $\phi:I\to\C$ be a parametrization of $\partial\Omega'$ for $I\subset\R$ an interval, and let
	$$
	\mathrm{dbl}=\{z\in\p \Omega': \exists s\neq t \;\textrm{such\,\,that\,}\; \phi(s)=\phi(t)=z\}.
	$$  
	Let $a',b'\in \partial\Omega'\setminus \mathrm{dbl}$ and $\phi: \ol{\bbH}\to\ol{\Omega'}$ 
	be a conformal map with $\phi(0)=a',\phi(\infty)=b'$.
	Let  $\sle'$ be an $\SLE_6$ on $(\Omega',a',b')$. We say that $(\Omega',a',b')$ satisfies  {\rm Property (S)} if    $\P[\dist(\sle',\mathrm{dbl})>0]=1$.
	If  $(\Omega',a',b')$ satisfies {\rm (S)} then Lemma~\ref{prop:coupling} holds for $(\Omega',a',b')$ with the same constants $c,C$.
\end{lemma}
\begin{proof}
	Suppose $\dist(\sle',\mathrm{dbl})>0$ a.s. Then $\P[\dist(\phi^{-1}(\sle'), \phi^{-1}(\mathrm{dbl}))<\delta]=o_\delta(1)$ for $\delta\in (0,1)$.
	Let $\bbH^\delta=\{z\in \bbH: \dist(z,\phi^{-1}(\mathrm{dbl}) ) >\delta \}$ and $\Omega^\delta=\phi(\bbH^\delta)$.   Then  $\P[\sle'\subset \Omega^\delta]=1-o_\delta(1)$.
	Since $\p \bbH^\delta$ is a simple curve, we see that $\Omega^\delta$ is a Jordan domain, thus  satisfying Lemma~\ref{prop:coupling}.  By the locality property of $\SLE_6$, the total variation distance between   the law of $\sle'$ and the $\SLE_6$ on $(\Omega^\delta,a',b')$ is $o_\delta(1)$. Since $c,C$ in Lemma~\ref{prop:coupling} are independent of $\delta$, letting $\delta\to 0$, we prove  Lemma~\ref{lem:coupling2}.
\end{proof}
In the notation of Lemma~\ref{lem:coupling2}, we say that $(\Omega',a',b')$ satisfies Property (W) if  $\sle'\cap \mathrm{dbl}=\emptyset$ a.s. The following lemma ensures that the complement of $\SLE_6$ hulls satisfies Property (W).
\begin{lemma}\label{lem:triple}
Suppose $\gamma$ is a chordal $\SLE_6$ as in Theorem~\ref{thm:conv}. Then a.s.\ there exists no point $p\in \gamma$ such that  $\gamma\setminus\{p\}$ is disconnected and $\gamma$ visits $p$ at least twice.
\end{lemma}
\begin{proof}
This is proved in  \cite[Remark 8.8]{Restriction}.
\end{proof}
Recall that $\cD_\bullet$ is a domain induced by a face with two ordered  marked points on its boundary.
By Lemma~\ref{lem:triple}, for $\bullet=i,j$, it is a.s.\ the case that $\cD_{\bullet}$ (after recentering at 0) satisfies Property (W). However, $\cD_\bullet$ does not satisfy Property      (S) because the two boundary marked points could be accumulation points of $\mathrm{dbl}$. (In fact, one can prove that the two boundary marked points a.s.\ are such accumulation points.) 
We overcome this issue by the following lemma.
\begin{lemma}\label{lem:local}
For $\alpha\in(0,1)$, let $\bbH_\alpha=(\bbH \setminus \alpha \D) \cap \alpha^{-1} \D $. 
Let $\sle^0$ and $\sle^\alpha$  be the chordal $\SLE_6$ on $(\bbH,0,\infty)$ and $(\bbH_\alpha, \alpha i, \alpha^{-1}i)$ respectively. Let $\bx\subset \H$ be a box. Let $\sigma$ and $\ol\sigma$ be the first and last, respectively, time that $\sle$ is contained in $\bx$. Define $\sigma_\alpha$ and $\ol\sigma_\alpha$ for $\sle^\alpha$ similarly. Then the total variation distance between $\sle|_{[\sigma,\ol\sigma]}$ and $\sle^\alpha|_{[\sigma_\alpha,\ol\sigma_\alpha]}$ as curves modulo monotone parametrizations is $o_\alpha(1)$ as $\alpha\to0$.
\end{lemma}
\begin{proof}
Let $\wt\sle^\alpha$  be a chordal $\SLE_6$ on $(\bbH\setminus \alpha\D, \alpha i, \infty)$. Let $\wt \sigma_\alpha$ and $\wt\tau^\alpha$  be the first and last, respectively, time that  $\wt\sle^\alpha$ is contained in $\bx$. 
We couple $\wt\sle^\alpha$ and $\sle^0$ such that when running them backward, the two curves agree until hitting $\sqrt\alpha\D$; this is possible by reversibility of SLE$_6$. With probability $1-o_\alpha(1)$, the remaining segments of the two curves will not touch $\bx$. 
Then  the total variation distance between  $\sle|_{[\sigma,\infty)}$ and $\wt\sle^\alpha|_{[\wt\sigma_\alpha,\infty)}$ as curves modulo monotone parametrizations is $o_\alpha(1)$. Similarly, the total variation distance between  $\wt\sle^\alpha|_{[\wt\sigma_\alpha,\wt\tau_\alpha]}$ and $\sle^\alpha|_{[\sigma_\alpha,\ol\sigma_\alpha]}$ as curves modulo monotone parametrizations is $o_\alpha(1)$.  This concludes the proof.
\end{proof}
 Now we are ready to prove \eqref{eq:coup}.
Let $\phi$ be defined as in Lemma~\ref{lem:coupling2} with $\cD_\bullet$ and $x_\bullet$ in place of $(\Omega',a',b')$  and $0$. Recall the notation in Lemma~\ref{lem:local}. 
We can define the analog of $x_\bullet, y_\bullet$ with $\sle^\bullet$ replaced by the  $\SLE_6$ on $( \phi (\bbH_\alpha), \phi(\alpha i ), \phi (\alpha^{-1}i) )$ and denote these two quantities by $x_\alpha,y_\alpha$. 
By Lemma~\ref{lem:local},  the total variation distance between the laws of $(x_\bullet,y_\bullet)$ and $(x_\alpha,y_\alpha)$ is $o_\alpha(1)$. 
	On the other hand, $( \phi (\bbH_\alpha), \phi(\alpha i ), \phi (\alpha^{-1}i) )$  satisfies the stronger property (S) rather than just (W) because the boundary is simple and smooth near $\phi(\alpha i ) $ and $ \phi (\alpha^{-1}i)$.
	Since $c,C$ in Lemma~\ref{lem:coupling2} are independent of $\alpha$, letting $\alpha\to 0$, we arrive at
	\begin{equation*}\label{e.before}
		\Eb{x_\bullet - \nor y_\bullet \bigm| \cD_{\bullet}} \le  C\left(\frac {2\eps}  {\intd}\right)^{c} \cdot  \eps^{7/4} \cdot \alpha_2 (2\eps,\zeta)= o_{\eps/\intd} (1) \frac{\eps^2}{\intd^{1/4}}.
	\end{equation*}This concludes the proof of \eqref{eq:coup} and hence of Proposition~\ref{l.XY}.

\section{Minkowski content for percolation pivotal points}\label{sec:piv}
This section is devoted to proving Proposition~\ref{prop:occu-piv} and Theorem~\ref{thm:piv}. 

Recall that $\SLE_\kappa(\rho)$ and $\SLE_\kappa(\rho_1;\rho_2)$ processes are variants of $\SLE_\kappa$ whose driving functions have forcing terms prescribed by \notion{force points} with certain \notion{weights}. $\SLE_\kappa(\rho)$ has a single force point of weight $\rho$, while  $\SLE_\kappa(\rho_1;\rho_2)$ has two force points of weight $\rho_1$ and $\rho_2$, respectively. We will not give the formal definition of these processes and refer instead to \cite[Section 2.2]{IG-I}, because we only use a few well-established facts about the processes developed in the framework of imaginary geometry \cite{IG-I}. 

We also recall the Brownian excursion on $\bbH$ from $0$ to $\infty$. See \cite[Chapter 2]{Lawler-Book} for the precise definition. By the theory of conformal restriction \cite{Restriction} 
 the left and right boundary of the Brownian excursion and those of $\SLE_6(2;2)$ have the same law.  (In fact,  the hull of both Brownian excursion and $\SLE_6(2,2)$ are the unique chordal restriction measure with exponent 1.)
Let $\cC$ denote the intersection of the left and right boundaries of the Brownian excursion, i.e.\ the set of cut points. By \cite[Theorem 4.7]{cutpoint}, in the notation of Definition~\ref{def:occu}, we have the following.
\begin{lemma}\label{lem:cut}
The occupation measure of $\cC$ a.s.\ exists and is defined by its $3/4$-dimensional Minkowski content.
\end{lemma}

\begin{remark}\label{rmk:cutpoint}
In \cite[Theorem 4.7]{cutpoint}, the notion of a cut point is defined via cut times. Namely, given a Brownian excursion $(E(t))_{t\geq 0}$ on $\bbH$ from 0 to $\infty$, 
the set of cut points of $E$ is defined by $\cC'=\{ \eta(t)\,:\, t\geq 0 \op{\,\,and\,\,}E((0,t))\cap E((t,\infty))=\emptyset\}$. 
However, it can be checked that $\cC'=\cC$ a.s.  The direction $\cC'\subset \cC$ a.s.\ is trivial. For the other direction, 
let $\cH$ be the hull of $E$, which has the same law as the hull of $\SLE_6(2,2)$. By the SLE duality, the interior of $\cH$ is a countable collection of simply connected open sets, ordered by the order in which they are first visited if we go from 0 to $\infty$ inside $\cH$. In particular, for each $p\in \cC$, $\cC\setminus\{p\}$ has  two components, one bounded and one unbounded, such that all the sets in the bounded component are ordered before the sets in the unbounded component. 
Suppose there exists 
$p\in \cC\setminus\cC'$.
Let $\cC_1$ and $\cC_2$ be the bounded and unbounded component of  $\cC\setminus\{p \}$, respectively. 
Let $t_1$ be the first time $E(t)=p$ and $q\in E((0,t_1))\cap E((t_1,\infty))$. 
Then $q\in \cC_1$ and there exists $t'>t_1$ such that $E(t')=q$.
For a rational $s\in (t_1,t')$, $E(s)$ a.s.\ is contained in a component $B$ of the interior of $\cH$.
Moreover, the closure of $B$ only has two points in $\cC$, one of which must be visited by $E$ twice.
Since there are only countably many such points, this can be ruled out by the strong Markov property of $E$ and the fact that a planar  Brownian motion a.s.\ does not visit any fixed point.  This gives $\cC\subset \cC'$ a.s.
\end{remark}

\smallskip

\begin{figure}
	\centering
	\includegraphics[scale=1]{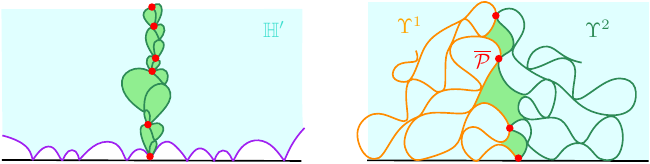}
	\caption{Left: The green curve is an SLE$_6(1,1)$ in $\H'$ (which is the domain in light blue) and has the law of an SLE$_6(2,2)$ viewed as a curve in $\H$. The points of intersection of its left and right boundaries (red) have the law of the cut points $\cC$ of a Brownian excursion in $\H$. Right: The region (light green) between the right boundary of $\Upsilon^1$ and the left boundary of $\Upsilon^2$ has the law of the region enclosed by an  SLE$_6(1,1)$ in $\H$.}\label{figure2}
\end{figure}

Run an $\SLE_{8/3}(2; -4/3)$ on $(\bbH,0,\infty)$ where the force points are at $0^-$ and $0^+$. 
Conditioning on this curve, run an $\SLE_{8/3}(-4/3;4/3)$ on the domain to its left.  Let $\bbH'\subset \bbH $ be the domain between these two curves. Conditioning on $\bbH'$, we run an $\SLE_6(1;1)$ on $(\bbH',0,\infty)$. Then by the rule of interacting flow lines in  \cite{IG-I},  the marginal law of this curve is an $\SLE_6(2;2)$  on $(\bbH,0,\infty)$ with force points at $0^+$ and $0^-$. See the left part of Figure \ref{figure2} for an illustration.

Let  $\Upsilon^1$ be an $\SLE_6(2)$ on $(\bbH,0,\infty)$ where the single force point is at $0^+$. Then $\Upsilon^1\cap \R_{>0} =\emptyset$.
Now conditioning on $\Upsilon^1$, let $\Upsilon^2$ be a chordal $\SLE_6$ from $0$ to $\infty$ on the domain to the right of $\Upsilon^1$.  By $\SLE$ duality (see \cite[Theorem~5.1]{Zhan-Duality} and \cite[Theorem~1.4]{IG-I})   the right boundary of $\Upsilon^1$ and the left boundary of $\Upsilon^2$ have the same joint law as the left and right boundary of an $\SLE_6(1,1)$. See the right part of Figure \ref{figure2} for an illustration. Denote their intersection by $\ol \cP$. Combined with the paragraph above, we have the following.
\begin{lemma}\label{lem:coupling}
There is a coupling of  $\cC$, $\bbH'$, and $\ol\cP$ such that $\ol\cP$  is independent of $\bbH'$, and $\cC$ is the image of $\ol\cP$ under a conformal map from $\bbH$ to $\bbH'$ fixing $0$ and $\infty$.
\end{lemma}
The next lemma links $\ol\cP$ to the set $\cP$ in Proposition~\ref{prop:occu-piv} and Theorem~\ref{thm:piv}.
\begin{lemma}\label{lem:P}
Recall $(\Omega,a,b,c,d)$ and $\sle^1,\sle^2$ in Theorem~\ref{thm:length2}. Let $\ol\sle^2$ be the time-reversal of $\sle^2$.
 There exist random times $\sigma$ and $\tau$ for $\Upsilon^1$ satisfying $0<\sigma<\tau<\infty$ with positive probability, such that the following hold on the event that $0<\sigma<\tau<\infty$.
\begin{enumerate}
\item \label{item:domain} The unbounded component of  $\bbH\setminus(\Upsilon^1([0,\sigma])\cup \Upsilon^1([\tau,\infty])$ can be conformally mapped to $\Omega$ with $(\Upsilon^1(\sigma), 0,\infty,  \Upsilon^1(\tau))$ mapped to $(a,b,c,d)$.
\item \label{item:resample} Conditioning on the realization of $\bbH\setminus(\Upsilon^1([0,\sigma])\cup \Upsilon^1([\tau,\infty])$, the joint law of the conformal image of $\Upsilon^1([\sigma,\tau])$ and  $\Upsilon^2$ is the same as the conditional  law of $(\sle^1,\ol\sle^2)$ conditioning on $E^c$. 
\end{enumerate}
\end{lemma}
\begin{figure}
	\centering
	\includegraphics[scale=1]{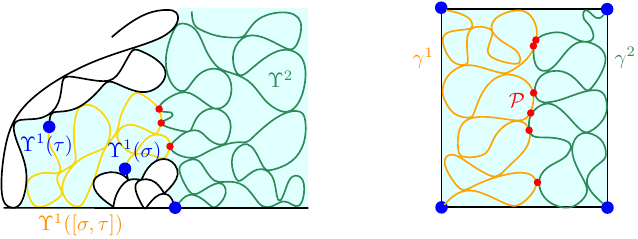}
	\caption{Illustration of Lemma \ref{lem:P}.}
\end{figure}
\begin{proof}
Let $t$ be the last time where $\op{Im} (\Upsilon^1)=1$. 
By \cite[Lemma 2.2]{Miller-Wu},  $\Upsilon^1[0,t]$ stays close to any deterministic smooth curve in $\bbH$ from the origin to a point on $\{z:\op{Im}z=1\}$ with positive probability.
Therefore, with positive probability  $\Upsilon^1$ reaches a time $s<t$ when 
Condition~\ref{item:domain} in Lemma~\ref{lem:P} is achieved with $s,t$ in place of $\sigma,\tau$. On the event that there is such a time $s\in(0,t)$, let $\sigma$ be the infimum of such times. Note that $\sigma>0$ a.s.\ since as $s\to 0$ the extremal distance between the arc from  $\Upsilon^1(s)$ to $0$ and the arc from $\infty$ to $\Upsilon^1(\tau)$ in the unbounded component of $\bbH\setminus(\Upsilon^1([0,s])\cup \Upsilon^1([\tau,\infty]))$ goes to zero. We set $\sigma=\infty$ if this event does not occur and  let $\tau=t\vee \sigma$.

Viewing $\Upsilon^{1},\Upsilon^{2}$ as two counterflow lines of different angles in the same  imaginary geometry (\cite{IG-I}), for $i=1,2$, conditioning on $\Upsilon^{i}$, the law of $\Upsilon^{3-i}$ is a chordal $\SLE_6$. Therefore the domain Markov property and reversibility of $\SLE_6$ yield that the same resampling property holds for  the conformal image of $\Upsilon^1([\sigma,\tau])$ and  $\Upsilon^2$.  By convergence of $(\sle^1_\mesh,\sle^2_\mesh)$ to $(\sle^1,\sle^2)$, the same resampling property holds for  $(\gamma^1,\ol \gamma^2)$ conditioning on $E^c$.  As explained in \cite[Appendix A]{MSW-meas}, this resampling property uniquely determines the law of the pair of curves.  Thus we conclude the proof.
\end{proof} 
Combining Lemmas~\ref{lem:cut}, \ref{lem:coupling}, and~\ref{lem:P}, 
we see that on the event $E^c$ the occupation measure of $\cP$ exists and is defined by its $3/4$-dimensional Minkowski content. The same argument works when conditioning on $E$. 
This gives Proposition~\ref{prop:occu-piv}.
\bigskip

The proof of Theorem~\ref{thm:piv} follows from the exact same argument as in the proof of Theorem~\ref{thm:equi-length}. We just need to replace one interface $\sle$ with the pair of  interfaces $\sle^1,\sle^2$. Here we only point out the substitutes of the ingredients in the argument in Section~\ref{sec:equiv}.  

Suppose we are in the coupling of Theorems~\ref{thm:length2} and~\ref{thm:embedding}.  Then Lemma~\ref{lem:mono2} and the 4-arm case of Lemma~\ref{lem:endpoint}  and  \eqref{eq:24'}  give  the analog of Lemmas~\ref{lem:12y}, \ref{lem:12pt}, and ~\ref{lem:two-arm}, in addition to \eqref{eq:beta}.
We can also adapt the concept of face in this setting, where the number of arcs becomes 4 instead of 2.
We use the same notations as in Section~\ref{subsec:couple}. Given $\bx$, let $\cA(\bx,\Omega)= \{\cP\cap \bx\neq \emptyset\}$. On the event $\cA(\bx,\Omega)$,  we trace $\gamma^1$ and $\gamma^2$ and their time-reversals from $a,b,c,d$ until first hitting $\bx$. This defines a face at $\Theta$ at $\bx$ induced by $(\gamma^1,\gamma^2)$. Moreover, $\cD(\bx,\Omega)$ and $\cG_\Omega(\bx,U)$ can be defined in the same way as in Section~\ref{subsec:couple}.   Then Lemmas~\ref{lem:GG} and \ref{lem:GG2} still hold with $k=2,3$ replaced by $k=4,5$. Now if we carry out the argument in Section~\ref{sec:equiv}, Theorem~\ref{thm:piv} will be reduced to the analog of \eqref{eq:coup}, which can still be proved by a coupling argument as in Lemma~\ref{prop:coupling} combined with approximation arguments as in Lemmas~\ref{lem:coupling2}, \ref{lem:triple}, and~\ref{lem:local}.

\section{Convergence of the percolation interface under natural parametrization}\label{sec:conv}
In this section we prove Theorem~\ref{thm:conv}. We pick {the constant $\con1$ so that Theorem \ref{thm:equi-length} holds}.
\begin{lemma}\label{lem:tight}
	The curve $\wh{\sle}_\mesh$ is tight for the $\rho$-metric. 
\end{lemma}
\begin{proof}
We proceed by contradiction and suppose $\wh\sle_\mesh$ is not tight. Then there exist $\delta_0>0$, $\mesh_n\downarrow 0$ and $\eps_n\downarrow 0$ such that
\begin{equation}\label{eq:osc}
\P[\osc(\eps_n; \wh{\sle}_{\mesh_n})>\delta_0]\ge \delta_0. 
\end{equation}
	where $\osc(\eps,f)=\sup_{|t-s|\le \eps} |f(t)-f(s)|$. Write $\wh \sle_{\mesh_n}$ as $\wh{\sle}_n$ and $\sle_{\mesh_n}$ as $\sle_n$ for simplicity.
	Given a realization of $\sle_n$, let $s_n,t_n$ be the smallest times such that  
	$|\wh{\sle}_n(s_n)-\wh \sle_n(t_n)|=\osc(\eps_n; \wh{\sle}_{n})$.
	Then the random variables $(s_n,t_n,\wh{\sle}_n(s_n),\wh \sle_n(t_n))$ are tight. By the Skorokhod embedding theorem, we can couple $\{\sle_n\}$ and $\sle$ such that (possibly) along a subsequence (which is still indexed by $n$ for the sake of simplicity), the  following occur a.s.
	\begin{enumerate}
		\item $\sle_n$ converge to $\sle$ in the topology of Theorem~\ref{thm:length};
		\item $s_n$ and $t_n$ converge to the same limit, which we denote by $\ft$;
		\item there are $x,y\in \ol\Omega$ such that $\wh \sle_n(s_n)\rightarrow x$ and $\wh \sle_n(t_n) \rightarrow y$.
	\end{enumerate}
	We first observe that both $x$ and $y$ are on the trace of $\sle$ a.s.  In fact, consider an open ball $B(z,r)$ where $z,r$ are both rational. (We call such balls \notion{rational balls}.) Then on the event that $x\in B(z,r)$, it holds a.s.\ that for sufficiently large $n$, $\sle_n\cap B(z, r)\neq\emptyset$. Given condition 1 in the coupling above, $\sle\cap B(z,r)\neq \emptyset$ a.s.\ (see Lemma~\ref{lem:stop}).
	
	By \eqref{eq:osc},  in the coupling above, $\P[|x-y|\ge \delta_0] \ge \delta_0$. Therefore we can find rational balls $B(z_1,r_1)$ and $B(z_2,r_2)$  such that the following event occurs with positive probability:
	\begin{enumerate}
		\item $x\in B(z_1,r_1)$ and $y\in B(z_2,r_2)$;
		\item $\max\{r_1,r_2\}<0.1 \delta_0$ and
		$B(z_1,2r_1)\cap B(z_2, r_2)=\emptyset$.
	\end{enumerate}
	We work on this event hereafter. Let $\sle_n$ and $\sle$  be parametrized as in Lemma \ref{lem:stop} and let
	\begin{align*}
		\rho^1_n&=\inf \{ t: \sle(t) \in B(z_1,r_1)\},\\
		\sigma^1_n&=\inf\{t>\rho^1_n: \sle_n(t) \notin B(z_1,2r_1) \},\\
		\lambda^1_n&=\inf\{t>\sigma^1_n: \sle_n(t) \in B(z_2,r_2)\}.
	\end{align*}
	Define $\rho^1, \sigma^1,\lambda^1$ similarly for $\sle$.  Note that $K_{\lambda^1}\setminus K_{\sigma^1}$ has non-empty interior a.s., where $K_\cdot$ is the hull process of $\sle$. Therefore  there exists a rational ball $B(z_3,r_3)\subset K_{\lambda^1}\setminus K_{\sigma^1}$ such that
	$$\sle\cap B(z_3,r_3)= \sle([\sigma^1,\lambda^1]) \cap B(z_3,r_3)\neq \emptyset.$$ 
	By Lemma~\ref{lem:stop}, we must have that  for all sufficiently large $n$, 
	$$\sle_n\cap B(z_3,r_3)= \sle_n([\sigma^1_n,\lambda^1_n]) \cap B(z_3,r_3)\neq\emptyset.$$
 Since $\om(\p B(z_3,r_3))=0$,	by Theorem~\ref{thm:equi-length}, $\lambda^1_n-\sigma^1_n\ge \tau_n( B(z_3,r_3)) \to c_*\om (B(z_3,r_3))>0$.   Since $\sle_n([0,\rho^1_n])  \subset \wh{\sle}_n([0,s_n])$ and $\wh \sle_n([0,\lambda_n^1])  \subset \sle_n([0,t_n])$, while $t_n-s_n\to 0$, we must have 	 $\sle_n([0,\sigma^1_n])  \subset \wh{\sle}_n([0,s_n])$.

	Now for  $k\ge 2$ we  define $\rho^k_n$ and $\sigma^k_n$ inductively by
	\begin{align*}
		\rho^k_n&=\inf \{ t>\sigma^{k-1}_n : \sle_n(t) \in B(z_1,r_1)\},\\
		\sigma^k_n&=\inf\{t>\rho^k_n: \sle_n(t) \notin B(z_1,2r_1) \}.
	\end{align*}
	Then $\rho^2_n<\infty$. By the same argument as above, for any fixed  $k$,  for $n$ sufficiently large,  we have  $\sle_n([0,\sigma^k_n])  \subset \wh{\sle}_n([0,s_n])$ hence $\rho^{k+1}_n<\infty$.
	This yields that with positive probability, $\sle$ returns to $B(z_1,r_1)$ after hitting $\p B(z_1,2r_1)$ infinitely many times.  {This contradicts the fact that  almost surely  $\sle$ is a continuous curve and never returns to $B(z_1,2r_1)$ after a certain finite time.}
\end{proof}
Before proving Theorem~\ref{thm:conv} we prepare two lemmas.
\begin{lemma}\label{lem:neighbor}
	Let $\sle, \wh \sle$ be as in Theorem~\ref{thm:conv} and recall that $\om$ denotes the occupation measure of $\sle$. For each fixed  $t>0$, on the event $\om(\Omega)>t$ let $K_t$ be the hull of $\wh \sle([0,t])$ (see Lemma~\ref{lem:mono}). Then $\om(\partial K_t)=0 $ a.s.
\end{lemma}
\begin{proof}
By Assertion~\ref{item:regular} above Definition~\ref{def:natural},  with probability 1,  $\om(\wh{\sle}(t))=0$ for all $t\ge 0$.  To prove Lemma~\ref{lem:neighbor},  it suffices to show that for any fixed square $Q\subset \Omega$ and $\delta>0$, we have 
$$\om (U^\delta)=0\quad \textrm{a.s., where}\quad U^\delta =Q\cap \partial K_t\setminus B(\wh{\sle}(t),\delta)$$  
For $0<r<0.1\eps< 0.01\delta$, let $U^{\delta,\eps}$ be the union of all  boxes  of  side length $\eps$ on $\eps\Z^2$   that has nonempty intersection with $U^\delta\neq \emptyset$.
 Let $U^{\delta,\eps}_r$ be the $r$-neighborhood of $U^{\delta,\eps}\cap \sle$.
Suppose we are in the coupling of Theorem~\ref{thm:embedding} so that $\sle$ is coupled with $\omega\in\cH_\Omega$. 
By Lemma~\ref{lem:mono}, for each $x\in\Omega$, if $x\in U^{\delta,\eps}_r$, then we have
	\begin{enumerate}
		\item the 2-arm event occurs for the annulus of $\cB(x,\eps)\setminus \cB(x,r)$;
		\item the 3-arm event occurs for the annulus of $\cB(x,0.5\delta)\setminus \cB(x,2\eps)$.
	\end{enumerate}
	By \eqref{eq:24'} and Lemma~\ref{lem:independent}, there exists a constant  $C=C(\delta)$ such that for any $x\in\Omega$,
	\[
	\P[x\in U^{\delta,\eps}_r]\le C(\delta) (r/\eps)^{1/4}\eps^{2/3}=C(\delta)r^{1/4} \eps^{5/12}.
	\] 
By the definition of $\om$ in Definition~\ref{def:occu},  almost surely $\om (U^{\delta,\eps})=\Mink_{7/4}(U^{\delta,\eps} \cap \sle)$. By Fatou's lemma and Definition~\ref{def:Mink}, 
\[
\E[\om(U^\delta)]\le\E[\om (U^{\delta,\eps})]	\le  \liminf_{r\to0} \E [r^{-1/4} \mathrm{Area}( U^{\delta,\eps}_r  )  ]\le C(\delta)\eps^{5/12}.
\] 
Sending $\eps\to 0$ we have $\E[\om(U^\delta)]=0$ and we are done.
\end{proof}
\begin{lemma}\label{lem:uniform}
In the setting of Theorem~\ref{thm:conv}, conditioning on $\gamma$, sample $U$  from $(0,\tau(\Omega))$ uniformly and sample $z\in\Omega$ according to $\tau(\cdot)/\tau(\Omega)$. 
Then $\wh \gamma(U)$ has the same law as $z$. 
\end{lemma}
\begin{proof}
For $t>0$, let $K_t$ be the hull of $\wh\gamma([0,t])$. Namely $A_t$ is the set $A$ in Lemma~\ref{lem:neighbor}. 
For $s\in [0,t)$, let $K_{s,t}=K_t\setminus K_s$. It is  proved in  \cite[Equation (27)]{Lawler-Minkowski}  that a.s.
\begin{equation}\label{eq:occu=Mink}
\om(\wh \sle[s,t])=t-s\qquad \textrm{for all } t>s\ge 0.
\end{equation}
By Theorem~\ref{thm:equi-length}, Lemma~\ref{lem:neighbor}, and \eqref{eq:occu=Mink}, for a fixed pair of $s,t$ we have $\tau(K_{s,t})=t-s$ a.s. For each $n\in \N$,  we can couple $z$ and $U$ such that both $\wh \gamma(U)$ and $z$ fall in $K_{i/n,(i+1)/n}$ for some $i<n$. By the continuity of $\wh \gamma$, we conclude the proof.
\end{proof}
\begin{proof}[Proof of Theorem \ref{thm:conv}]
	According to Lemma~\ref{lem:tight}, by possibly restricting to a subsequence, we can assume that $\wh{\sle}_\mesh$ converge a.s.\ to a curve $\sle'$ in the $\rho$-metric. Furthermore, we can assume that along this subsequence $(\tau_\mesh,\sle_\mesh)$ converges to $(\tau,\sle)$ a.s.\  in the topology of Theorem~\ref{thm:length}, where $\sle$ is a chordal $\SLE_6$ (viewed as a curve modulo reparametrization of time). In this coupling, $\sle'$ is a parametrization of $\sle$ with total length $\tau(\Omega)$. Let $\wh{\sle}$ be $\sle$ with its natural parametrization. It suffices to show that $\sle'=\wh\sle$ a.s.
	
	Conditioning on all the data above, we can sample a random time $U$  uniformly in $(0,\tau(\Omega))$ and a random edge $e_\mesh$ according to $\tau_\mesh$ on $\sle_\mesh$  such that $\P[\wh{\sle}_\mesh(U)\in e_\mesh] \ge 1-o_\mesh(1)$.  Notice that $\wh{\sle}_\mesh(U)$ converges to $\sle'(U)$ a.s. Therefore $e_\mesh$  converge to $\sle'(U)$ a.s. Here we  identify edges on the hexagonal lattice  with their midpoints since the difference is negligible in the scaling limit.  On the other hand, since $\tau_\mesh$  converge to $\tau$ in the Prokhorov metric a.s., by Lemma~\ref{lem:uniform} we have that $e_\mesh$  converges in law to $\wh{\sle}(U)$.  This implies that $\sle'(U)$ and $\wh \sle(U)$ are equal in law. 
	
	Given any fixed $t>0$, on the event that $t\in (0,\tau(\Omega))$, let $A$ be the hull of $\wh \sle([0,t])$ and $A_\eps=\{z\in \C: \dist(z,A)) \le \eps \}$.
	Since $\tau(A\setminus \partial A ) \le t$ and (by Theorem~\ref{thm:equi-length} and Lemma~\ref{lem:neighbor}) $\tau(\partial A) =0$, we have $\lim_{\eps\to 0}  \tau (A_\eps)=t$. 
	Thus for all $\delta>0$ with $t+2\delta<\tau(\Omega)$, we can find a (random) $\eps>0$ a.s.\  such that 
	$\tau(A_\eps)\le t+\delta$. Therefore $\limsup_{\mesh\to 0} \tau_\mesh(A_\eps) \le t+\delta$. Hence for sufficiently small $\mesh$ there exists $t^\delta_\mesh \in [t,t+2\delta]$ such that  $\wh\sle_\mesh(t^\delta_\mesh) \notin A_\eps$. By possibly passing to a subsequence, we can assume $\lim_{\mesh\to 0}t^\delta_\mesh\to t^\delta$. Sending $\mesh\to0$, we have $\sle'(t^\delta)\notin \wh \sle([0,t])$. Since $\sle'$ and $\wh \sle$ have the same range, we must have $\wh \sle(t)\in \sle'([0,t^\delta])$.   Letting $\delta\to 0$, we have $\wh\sle(t)\in \sle'([0,t])$. By considering rational $t$'s and then using the continuity of $\wh \sle$ and $\sle'$, we see that a.s.\ $\wh\sle([0,t])\subset \sle'([0,t])$ for all $t\in (0,\tau(\Omega))$. 
	Combined with $\sle'(U)\overset{d}{=}\wh \sle(U)$, we have $\sle'(U)=\wh \sle(U)$  hence $\sle'=\wh\sle$ a.s.
\end{proof}

\end{document}